\documentclass[12pt,reqno]{amsart}
\usepackage[activeacute,english]{babel}

\usepackage[utf8]{inputenc}

\usepackage[english]{babel}
\usepackage{amsfonts,amsmath,amsthm,bm,amssymb}
\usepackage{graphics,tikz,float}

\usepackage{hyperref,cleveref}
\usepackage{autonum} 

\DeclareRobustCommand{\SkipTocEntry}[5]{}

\usepackage[margin=0.9in]{geometry}
\newcommand{\PP}{{\mathcal P}}

\newcommand{\N}{{\mathbb N}}
\newcommand{\R}{{\mathbb R}}
\newcommand{\Z}{{\mathbb Z}}

\newcommand{\Sph}{{\mathbb S}}

\newcommand{\eper}{\operatorname{per}}

\newcommand{\intp}{\int_{-\pi}^{\pi}}
\newcommand{\essinf}{\operatorname{ess}\inf}
\newcommand{\esssup}{\operatorname{ess}\sup}

\newtheorem{theorem}{Theorem}[section]
\newtheorem{lemma}[theorem]{Lemma}

\newtheorem{remark}[theorem]{Remark}

\numberwithin{equation}{section}

\allowdisplaybreaks[1]


\begin{document}

\title[Strict rearrangement inequalities]
{Strict rearrangement inequalities: \\ 
 nonexpansivity and periodic Gagliardo seminorms}

\author[G. Csat\'o]{Gyula Csat\'o}
\address{G. Csat\'o \textsuperscript{1,2}
\newline
\textsuperscript{1}
Departament de Matem\`{a}tiques i Inform\`{a}tica,
Universitat de Barcelona,
Gran Via 585,
08007 Barcelona, Spain
\newline
\textsuperscript{2}
Centre de Recerca Matem\`atica, Edifici C, Campus Bellaterra, 08193 Bellaterra, Spain.}
\email{gyula.csato@ub.edu}

\author[A. Mas]{Albert Mas}
\address{A. Mas \textsuperscript{1,2}
\newline
\textsuperscript{1}
Departament de Matem\`atiques,
Universitat Polit\`ecnica de Catalunya,
Campus Diagonal Bes\`os, Edifici A (EEBE), Av. Eduard Maristany 16, 08019
Barcelona, Spain
\newline
\textsuperscript{2}
Centre de Recerca Matem\`atica, Edifici C, Campus Bellaterra, 08193 Bellaterra, Spain.}
\email{albert.mas.blesa@upc.edu}

\date{\today}
\thanks{The two authors are supported by the Spanish grants PID2021-123903NB-I00 and RED2022-134784-T funded by MCIN/AEI/10.13039/501100011033 and by ERDF ``A way of making Europe'', and by the Catalan grant 
2021-SGR-00087. The first author is in addition supported by the Spanish grant PID2021-125021NA-I00.
This work is supported by the Spanish State Research Agency, through the Severo Ochoa and Mar\'ia de Maeztu Program for Centers and Units of
Excellence in R\&D (CEX2020-001084-M)}

\begin{abstract}
This paper deals with the behavior of the periodic Gagliardo seminorm under two types of rearrangements, namely under  a periodic, and respectively a cylindrical, symmetric decreasing rearrangement. Our two main results are P\'olya-Szeg\H o type inequalities for these rearrangements. We also deal with the cases of equality.

Our method uses, among others, some classical nonexpansivity results for rearrangements for which we provide some slight improvements. Our proof is based on the ideas of [Frank and Seiringer, {\em Non-linear ground state representations and sharp Hardy inequalities}, J. Funct. Anal., 2008], where a new proof to deal with the cases of equality in the nonexpansivity theorem was given, albeit in a special case involving the rearrangement of only one function. 

\end{abstract}

\maketitle

\tableofcontents


\section{Introduction}

\subsection{The periodic Gagliardo seminorm and rearrangements}\label{ss rearrGagl}

The  goal of this paper is to establish P\'olya-Szeg\H o type inequalities for Gagliardo seminorms in a periodic setting. Let us start by defining, for $0<s<1$, $1\leq p<+\infty$, and a function $u:\R^n\to\R$ which is $2 \pi$-periodic in the variable $x_1$, the periodic Gagliardo seminorm
\begin{equation}\label{eq:def of periodic Gagl seminorm}
  [u]_{W^{s,p}}^{\eper}:=\left(\int_{\{x\in\R^n:\,-\pi<x_1<\pi\}}dx\int_{\R^n}dy\,\frac{|u(x)-u(y)|^p}{|x-y|^{n+s p}} \right)^{1/p};
\end{equation}
throughout this work we will use the notation 
$\R^n=\{x=(x_1,x')\in\R\times\R^{n-1}\}$.
This seminorm has been used in the literature in some specific situations. For instance, if $E\subset\R^n$ is $2 \pi$-periodic in the variable $x_1$ then, taking $u$ equal to the characteristic function~$\chi_E$ and $p=1$, we have that 
$\PP_s(E):=[\chi_E]_{W^{s,1}}^{\eper}$ is the periodic fractional perimeter introduced in 
\cite{Davila}. $\PP_s(E)$ was successfully used in \cite{CCM Delaunay} to construct periodic surfaces with constant nonlocal mean curvature. Whereas if $p=2$ and $n=1$ then 
$[\,\cdot\,]_{W^{s,1}}^{\eper}$ coincides with the kinetic part of the Lagrangian (or energy) studied in \cite{CCM Semilinear Var} to address the variational formulation, symmetry, and regularity properties of periodic solutions to semilinear equations involving the fractional Laplacian.

The celebrated P\'olya-Szeg\H o inequality states that, for sufficiently smooth open sets $\Omega\subset\R^n$,
\begin{equation}
 \label{eq:gy:intro:Polya Szego classic}
  \|\nabla (u^{*\,n})\|_{L^p(\Omega^{*\,n})}\leq \|\nabla u\|_{L^p(\Omega)}\quad\text{ for all }u\in W_0^{1,p}(\Omega),
\end{equation}
where $u^{*\,n}$ denotes the radially symmetric decreasing (Schwarz) rearrangement of the function $u:\Omega\subset\R^n\to \R$ and $\Omega^{*\,n}$ denotes the ball centered at the origin and with the same volume as $\Omega$. The vanishing boundary condition is important here.  Without it the inequality does not hold true in general; see the work of Kawohl \cite[Section II.4, Example 2.2]{Kawohl book}. We are interested in analogues of this inequality for the periodic Gagliardo seminorm.
We will establish two theorems dealing with its behavior under two different types of rearrangements. Let us briefly introduce these rearrangements; we refer to \Cref{section preliminaries} for the precise definition. The first one is the \textit{periodic rearrangement} $u^{*\eper},$ which is defined as follows: first perform a Steiner symmetrization with respect to the hyperplane $\{x_1=0\}$ of the function $u$ resticted to $(-\pi,\pi)\times\R^{n-1}$, and then extend this function to $\R^n$ in a $2\pi$-periodic way with respect to the variable $x_1$. The second rearrangement is the \textit{cylindrical rearrangement} with respect to the $x_1$-axis. It is denoted by $u^{{*\,n,1}}$ and it is obtained by performing the Schwarz rearrangement in $\R^{n-1}$ of $u(x_1,\cdot)$ for each frozen value of $x_1$.     

The behavior of the periodic Gagliardo seminorm under these two rearrangements has been addressed for the first time by D\'avila, del Pino, Dipierro, and Valdinoci \cite{Davila}, but only in the special case $p=1$ and if $u$ is a characteristic function. Their only result in this direction is \cite[Proposition 13]{Davila} which  establishes that $\PP_s(E^{{*\,n,1}})\leq \PP_s(E)$, 
where $E^{{*\,n,1}}:=\{(\chi_E)^{{*\,n,1}}>0\}$. We note that their proof also works for any function $u$ (and not only for $u=\chi_E$) as long as $p=1$ or $p=2,$ with slight modifications, however not for other values of $p.$ The same argument has also been used in \cite[Lemma 4.2]{Malchiodi Novaga Pagliardini}. Regarding the periodic rearrangement, the authors of \cite{Davila} only conjectured that $\PP_s(E^{*\eper})\leq \PP_s(E)$, where $E^{*\eper}:=\{(\chi_E)^{*\eper}>0\}$. We proved this conjecture in \cite{CCM Delaunay}, as well as the analogous inequality for the Gagliardo seminorm in the case $p=2$ and $n=1$ in \cite{CCM Semilinear Var}. Indeed, more general seminorms than \eqref{eq:def of periodic Gagl seminorm} are considered in this last work;  see  \Cref{remark:comparison n is one}. 

Let us briefly comment why these special cases are simpler, for both cylindrical and periodic rearrangements. First of all, for $u$ and $v$ general, the factorization 
 $$
   |u(x)-v(y)|^p=u^p(x)+v^p(y)-2u(x)v(y)
 $$
 only holds for $p=2$. We mention that it also holds if $p=1,$ $u=\chi_{E},$ and $v=\chi_{F}$ for any $E,\,F\subset\R^n.$ Hence in these cases one can directly apply  Riesz rearrangement inequalities to obtain P\'olya-Szeg\H o inequalities. Secondly, we will have to deal with sections of the kernel $|x-y|^{-(n+sp)}$ appearing in the seminorm \eqref{eq:def of periodic Gagl seminorm}, namely, the function 
$$t\in(0,+\infty)\mapsto g(t):=(t^2+a^2)^{-(n+ps)/2},$$ where $t=|x_1-y_1|$ and $a=|x'-y'|.$ If $n=1$ then $a=0$ and the function $g$ is convex. However, if $n>1$ then $g$ is concave for $t$ near the origin whenever $a\neq 0$, and this causes considerable difficulties when dealing with the periodic rearrangement. 

Our first main result is a full characterization of the behavior of the seminorm~$[\,\cdot\,]_{W^{s,1}}^{\eper}$ under periodic and cylindrical rearrangements. In particular, it includes the above mentioned partial results. In order to simplify the exposition, we only give here a brief summary of our main \Cref{thm:gy:periodic nonlocal Polya,thm:gy:Polya for n minus 1 Schwarz}, where we also deal with the cases of equality in the inequalities \eqref{eq:gy:intro:rearrangement general}.

\begin{theorem}
\label{thm:gy:intro:rearrangement general}
The periodic Gagliardo seminorm does not increase under periodic and cylindrical rearrangements. That is to say, if  $0<s<1$ and  $1\leq p<+\infty$, then
\begin{equation}
 \label{eq:gy:intro:rearrangement general}
  [u^{*\eper}]_{W^{s,p}}^{\eper}\leq [u]_{W^{s,p}}^{\eper}
  \quad\text{ and }\quad
  [u^{{*\,n,1}}]_{W^{s,p}}^{\eper}\leq [u]_{W^{s,p}}^{\eper}
\end{equation}
for all $u:\R^n\to\R$ which is measurable and $2 \pi$-periodic in the variable $x_1$.
\end{theorem} 

We emphasize that dealing with the cases of equality is considerably harder than establishing these inequalities; the reader will find the full description of the cases of equality in \Cref{subsection Gagliardo seminorms}.

Let us now make some comments on \Cref{thm:gy:intro:rearrangement general} in comparison with the known results in the nonperiodic local, nonperiodic nonlocal, and periodic local cases. We first concentrate on the nonperiodic setting.  The nonlocal counterpart of the classical P\'olya-Szeg\H o inequality \eqref{eq:gy:intro:Polya Szego classic} was first proven by Almgren and Lieb \cite{Almgren-Lieb} and states that
\begin{equation}
\label{eq:gy:intro:frac Polya Szego}
  [u^{*\,n}]_{W^{s,p}(\R^n)}\leq [u]_{W^{s,p}(\R^n)},\quad\text{ where }[u]_{W^{s,p}(\R^n)}^p:=\int_{\R^n}dx \int_{\R^n}dy\frac{|u(x)-u(y)|^p}{|x-y|^{n+s p}}.
\end{equation}
When studying the cases of equality, there is a major difference between the local and nonlocal case if $p>1.$ It is easy to see that equality in the classical P\'olya-Szeg\H o inequality \eqref{eq:gy:intro:Polya Szego classic} does not force $u$ to be equal to $u^{*\,n}$ modulo translations. 
Let us illustrate this with a simple example for $n=1$. Take a symmetric decreasing function $v$ such that $\operatorname{supp}(v)\subset (-2,2)$ and $v(x)=1$ for all $x\in [-1,1].$ Then, take $0<\epsilon<1$ and another nontrivial symmetric decreasing function $\varphi$ with $\operatorname{supp}(\varphi)\in (-\epsilon,\epsilon)$, and define $\varphi_t(x):=\varphi(x-t).$ Clearly, for every $1\leq p<+\infty$ and every $t\in (\epsilon-1,1-\epsilon),$ it holds that
\begin{equation}\label{eq:gy:intro:example Kawohl}
  \int_{-2}^2|(v+\varphi_t)'|^p
  =\int_{-2}^2|(v+\varphi_0)'|^p
  =\int_{-2}^2|((v+\varphi_t)^{*1})'|^p,
\end{equation}
but any translation of $v+\varphi_t$ differs from $(v+\varphi_t)^{*\,1}$ if $t\neq 0$.  However, if one excludes \textit{flat} (i.e., horizontal) parts of the graph of $u,$ then equality in the P\'olya-Szeg\H o inequality \eqref{eq:gy:intro:Polya Szego classic} does indeed force  $u$ to be equal to its rearrangement modulo translations. A first general theorem for a fairly big class of functions excluding \textit{flat} parts (for instance analytic functions) was given by Kawohl \cite{Kawohl ARMA}. This was generalized by Brothers and Ziemer \cite{Brothers Ziemer}, of which a simpler proof was later given by Ferone and Volpicelli \cite{Ferone Volpicelli}. 

In contrast to the local case, equality in the nonlocal P\'olya-Szeg\H o inequality 
\eqref{eq:gy:intro:frac Polya Szego} is sufficient to conclude that a translation of $u$ agrees with its rearrangement, if $p>1.$ Whereas for $p=1$ the conclusion in case of equality in \eqref{eq:gy:intro:frac Polya Szego} remains the same as in the local case: all superlevel sets of $u$ must be balls, but not necessarily centered at the same point (as in example \eqref{eq:gy:intro:example Kawohl}); see Frank and Seiringer  \cite[Appendix A]{Frank Seiringer}. Our characterization of the cases of equality in \Cref{thm:gy:intro:rearrangement general} is analogous, that is, if $p>1$ the function must be a translate of the corresponding rearranged function, whereas if $p=1$ the superlevel sets must be translates of the corresponding rearranged superlevel sets, but the translation might depend on the level.

When it comes to the rearrangement $u^{*\eper}$ of a function $u$ periodic in the variable $x_1$, there is no interest in studying the periodic version of the local P\'olya-Szeg\H o inequality \eqref{eq:gy:intro:Polya Szego classic} for $n=1$. Indeed, the inequality remains trivially true in each period. To see this, assume that $u$ is $2\pi$-periodic and with $u\in W^{1,p}(0,2\pi)$. Then simply consider the periods  
$\Omega_a:=(a+2k \pi,a+2(k+1)\pi),$ $k\in\Z,$ for some $a$ with $|u(a)|=\min |u|$. As $|u|-\min |u|\in W^{1,p}_0(\Omega_a),$ one can simply work separately in each period. A similar argument works for  $n=2$,  \textit{if} we require constant periodic boundary conditions $u(-\pi,x_2)=u(\pi,x_2)=c$ for all $x_2\in\R$, since the P\'olya-Szeg\H o inequality \eqref{eq:gy:intro:Polya Szego classic} remains  true also for Steiner symmetrization with vanishing boundary conditions; see Kawohl \cite[Section II.7, Corollary 2.32]{Kawohl book}. 

In this regard, we were not able to find in the literature a result that states or yields the following inequality: if $\Omega=(-\pi,\pi)\times\R,$  then 
\begin{equation}\label{eq:gy:intro:periodic local polya}
  \int_{\Omega}|\nabla (u^{*})|^p\leq \int_{\Omega}|\nabla u|^p
  \quad\text{for all $u\in W^{1,p}(\Omega)$ with  $u(-\pi,x_2)=u(\pi,x_2)$ for all $x_2\in\R,$}
\end{equation}
where $u^*$ is the Steiner rearrangement with respect to $\{x_1=0\}.$ This inequality can be proven using a periodic  polarization as introduced in Friedberg and Luttinger \cite{FriedLutt} and proceeding as in Brock and Solynin \cite[Sections 5, 6, and 8]{Brock Solynin}. 
The only explicit result of a periodic local P\'olya-Szeg\H o type inequality that we have found is the  following result by Cox-Kawohl \cite[Proposition 2.1]{Kawohl circular}: 
Let $\Omega=(-\pi,\pi)\times(0,1)$ and let $u$ satisfy 
$u(-\pi,x_2)=u(\pi,x_2)$ for all $x_2\in (0,1).$ Then,
$$
  \int_{\Omega}dx_1dx_2\,x_2\left(\frac{1}{x_2^2}(u^*)_{x_1}^2+(u^*)_{x_2}^2\right)
  \leq
  \int_{\Omega}dx_1dx_2\,x_2\left(\frac{1}{x_2^2}u_{x_1}^2+u_{x_2}^2\right),
$$
where $u_{x_1}$ and $u_{x_2}$ denote the partial derivatives of $u$.
 
The study of periodic versions  of the nonlocal P\' olya-Szeg\H o inequality \eqref{eq:gy:intro:frac Polya Szego}, in contrast to the local case,  becomes highly nontrivial even in dimension one. Firstly, the Gagliardo seminorm (or other nonlocal energies) feels the changes in all other periods and one cannot simply localize the problem to one period. Secondly, even if one could study only one period, there is no useful inequality of the type \eqref{eq:gy:intro:frac Polya Szego} for bounded domains $\Omega,$ i.e., when the double integral in the seminorm \eqref{eq:gy:intro:frac Polya Szego} is replaced by $\Omega\times\Omega.$  Indeed, symmetrization can increase the Gagliardo seminorm on domains; see for instance Li and Wang \cite{Li Wang}. For these two reasons one has to deal with the presence of competing terms when treating nonlocal periodic energy functionals. We have overcome these difficulties (in the case of the periodic rearrangement) by combining the following two ingredients via the Laplace transform: the nonexpansivity \Cref{theorem:gy:non expansivity under J} below and the monotonicity of the fundamental solution of the heat equation with periodic boundary conditions.


\subsection{A general nonexpansivity result for rearrangements}

Our proof of the P\'olya-Szeg\H o type inequalities in \Cref{thm:gy:intro:rearrangement general} strongly relies on the two very general nonexpansivity inequalities described in \Cref{theorem:gy:J in R n-1,theorem:gy:non expansivity under J} below. We present them here since they are of independent interest and contain some additions to earlier known versions. For instance, in there we also treat the case of equality for the convolution kernel equal to the absolute value $t\mapsto |t|$, which is not strictly convex. Moreover, our method to prove them is a generalization of the one used in  Frank and Seiringer
\cite[Lemma A.2]{Frank Seiringer}, who gave en elegant and quite elementary new proof of the case of equality in such a nonexpansivity inequality.

The classical nonexpansivity property of the $L^p$ distance of two functions under rearrangement can be seen as a generalization of the Riesz rearrangement inequality, which states
that
\begin{equation}
\label{eq:gy:intro:classical Riesz rearrangement}
  \int_{\R^n}dx\int_{\R^n}dy\,u(x)g(x-y)v(y)
  \leq 
  \int_{\R^n}dx\int_{\R^n}dy\,u^{*\,n}(x)g^{*\,n}(x-y)v^{*\,n}(y).
\end{equation}
Under the assumption that $g=g^{*\,n}$ is radially decreasing in the whole space, the cases of equality were first characterized by Lieb \cite{Lieb strict rearrang}.
Let us now state in the following theorem the generalization of this inequality, and we shall comment on its history and other versions thereafter. As stated here, it is the generalization of
\cite[Lemma A.2]{Frank Seiringer} to two functions. In the theorem, $|\Omega|$ denotes the measure of a set 
$\Omega\subset\R^n$.

\begin{theorem}
\label{theorem:gy:J in R n-1}
Let $J$ be a nonnegative, convex function in $\R$ with $J(0)=0$, and let $g\in L^1(\R^n)$ be a nonnegative function. Then,
\begin{equation}\label{cyl.rear.ineq} 
 \int_{\R^{n}}dx\int_{\R^{n}}dy\,J(u^{*\,n}(x)-v^{*\,n}(y))g^{*\,n}(x-y)
 \leq
 \int_{\R^{n}}dx\int_{\R^{n}}dy\,J(u(x)-v(y))g(x-y)
\end{equation}
for all pairs of nonnegative measurable functions $u,v:\R^n\to [0,+\infty)$ such that the right hand side of \eqref{cyl.rear.ineq}  is finite, and $|\{u>\tau\}|$ and $|\{v>\tau\}|$ are finite for all $\tau>0.$ 

If, in addition, $J$ is strictly convex and $g$ is a radially symmetric decreasing\footnote{Here, and throughout the whole paper, we call {\em decreasing} what some other authors call \em{striclty decreasing}.} function, then equality in \eqref{cyl.rear.ineq} holds if and  only if one of the following two cases occur: 
\begin{itemize}
\item[$(i)$] One of the two functions $u$ or $v$ is zero almost everywhere in $\R^n$ (and the other one can be anything).

\item[$(ii)$] There exists $a\in\R^{n}$ such that $u(x)=u^{*\,n}(x-a)$ and  $v(y)=v^{*\,n}(y-a)$  for almost every $x\in\R^{n}$  and almost every $y\in\R^{n}.$ 
\end{itemize}

If $J(t)=|t|$ and $g=g^{*\,n}$ is decreasing, then equality in \eqref{cyl.rear.ineq} holds if and only if for every  $\tau\in (0,\min\{\operatorname{ess}\sup u,\operatorname{ess}\sup v\})$ there exists $z_{\tau}\in\R^n$ such that 
$$
  \{u>\tau\}=\{u^{*\,n}>\tau\}-z_{\tau},\quad \{v>\tau\}=\{v^{*\,n}>\tau\}-z_{\tau}\quad\text{up to sets of measure zero.}
$$

\end{theorem}

The description of the cases of equality for strictly convex $J$ is a special case of Burchard and Hajaiej \cite[Theorem 2]{Burchard Hajaiej}, which dealt with more general convolution kernels than a convex function $J$ and multiple integrals\footnote{
To compare with \cite[Theorem 2]{Burchard Hajaiej}, we must mention that if $J$ is convex then $F(y_1,y_2):=-J(y_1-y_2)$ is supermodular.  See \cite{Brock Supermodular} for the case when $J$ is a $C^2$ function.}. 
However, \cite{Burchard Hajaiej} does not cover the cases of equality in the case $J(t)=|t|.$ Note also that \cite{Burchard Hajaiej} states the inequality only in the case that $g=g^{*\,n}$. To the best of our knowledge, the method of polarization used in \cite{Burchard Hajaiej} can only work in that special case when $g=g^{*\,n}$.

Note that we characterize the cases of equality only under an additional assumption on $g.$
For a general $J$ there is probably little hope of a simple characterization of the cases of equality  if we drop the assumption $g=g^{*\,n}.$ As far as we know, the only case where this  has been done is Burchard \cite{Burchard A} for $J(t)=|t|^2$ ---in which case \Cref{theorem:gy:J in R n-1}  reduces to the Riesz rearrangement inequality by expanding the product $(u(x)-v(y))^2=u(x)^2+v(y)^2-2u(x)v(y)$. In the next paragraph, when we talk about the cases of equality, we will always assume that $g=g^{*\,n}$ is radially decreasing.

Observe that for $g=\delta$ (the Dirac delta function understood as a measure) and $J(t)=|t|^p$ \Cref{theorem:gy:J in R n-1}
 is the well known nonexpansivity inequality for the $L^p$-norm, whereas for $J(t)=|t|^2$ it is the Riesz rearrangement inequality, as mentioned above. The inequality involving general $g$ and $J$ is due to Almgren and Lieb \cite[Corollary 2.3]{Almgren-Lieb} under the additional (but unnecessary) assumptions $J(-t)=J(t)$ and $J(u),\,J(v)\in L^1(\R^n).$  A proof of the inequality without these two assumptions and for more general $(k,n)$-Steiner rearrangements, but requiring $g=g^{*\,n}$,  can be found in  Brock and Solynin \cite[Lemma 8.2]{Brock Solynin}. This last proof uses polarization techniques and, therefore, only works  if $g=g^{*\,n}$. 
In the special case $u=v$ and $g=g^{*\,n}$ the two unnecessary assumptions were also removed by Frank and Seiringer \cite[Lemma A.2]{Frank Seiringer}, and indeed their proof also works for general $g$.
A characterization of the cases of equality, again in the special case $u=v,$ for strictly convex $J$, can be found in \cite[Theorem 8.1]{Brock Solynin}.
In addition, \cite{Frank Seiringer} characterizes the cases of equality under the assumption $u=v$ (not being aware of the work of \cite{Brock Solynin}), for strictly convex $J,$  but also for $J(t)=|t|$. We however need the inequality, and in particular the characterization of the equality cases, in the most general form when $u$ is not necessarily equal to $v.$ This generalization is necessary for our application to Gagliardo seminorms, namely, to the second inequality in \Cref{thm:gy:intro:rearrangement general}. In there we deal with cylindrical rearrangement and we will apply \Cref{theorem:gy:J in R n-1}
 in $\R^{n-1}$ to the functions $x'\mapsto u(x_1, x')$ and $v(y'):=u(y_1,y')$ for two different frozen values $x_1,\, y_1\in\R$.

For the proof of the first inequality in \Cref{thm:gy:intro:rearrangement general}, we will also need a version of the nonexpansivity \Cref{theorem:gy:J in R n-1}
 on the circle. This is the content of \Cref{theorem:gy:non expansivity under J}. There are some subtle differences with the nonperiodic result in the hypothesis, proof, and conclusion.  In its statement, $u^*$ denotes the Steiner rearrangement of $u$ with respect to the origin in the interval $(-\pi,\pi)$, which is the same as the Schwarz rearrangement in dimension one. Moreover, for a function $g$ which is $2\pi$-periodic in $\R$, recall that $g^{*\eper}$ denotes the $2\pi$-periodic extension of the function $(g\chi_{(-\pi,\pi)})^*$ (i.e., the Schwarz rearrangement of $g$ once restricted to one period); see \Cref{section preliminaries} for the detailed definition.

\begin{theorem}
\label{theorem:gy:non expansivity under J}
Let $J$ be a nonnegative, convex function in $\R$, and let $g$ be a nonnegative measurable $2\pi$-periodic function with $g\in L^1(-\pi,\pi)$. Then,
\begin{equation}\label{rear:Euv}
  \int_{-\pi}^{\pi}dx\int_{-\pi}^{\pi}dy\,J(u^*(x)-v^*(y))g^{*\eper}(x-y)
  \leq
  \int_{-\pi}^{\pi}dx\int_{-\pi}^{\pi}dy\,J(u(x)-v(y))g(x-y)
\end{equation}
for every pair of nonnegative measurable functions $u,v:(-\pi,\pi)\to\R$ such that the right hand side of \eqref{rear:Euv} is finite. 

If, in addition, $J$ is strictly convex, $g=g^{*\eper}$, and $g$ is decreasing in $(0,\pi)$,  then equality holds in \eqref{rear:Euv} if and only if one of the two cases occur:
\begin{itemize}
\item[$(i)$] One of the two functions $u$ or $v$ is constant almost everywhere in $(-\pi,\pi)$ (and the other one can be anything).

\item[$(ii)$] If $u$ and $v$ are extended to $\R$ in a $2\pi$-periodic way, then there exists $z\in\R$ such that $u(x)=u^{*\eper}(x+z)$ and $v(y)=v^{*\eper}(y+z)$ for almost every $x,\,y\in\R$.
\end{itemize}

If $J(t)=|t|$, $g=g^{*\eper}$, $g$ is decreasing in $(0,\pi),$ and $u$ and $v$ are extended to $\R$ in a $2\pi$-periodic way, then equality in \eqref{rear:Euv} holds if and only if for every $\tau\in (\operatorname{ess}\inf u,\operatorname{ess}\sup u)\cap (\operatorname{ess}\inf v,\operatorname{ess}\sup v)$ there exists $z_{\tau}\in\R$ such that 
\begin{equation}\label{rear:|t|.sets}
  \{u>\tau\}=\{u^{*\eper}>\tau\}-z_{\tau},\quad 
  \{v>\tau\}=\{v^{*\eper}>\tau\}-z_{\tau}\quad\text{up to sets of measure zero.}
\end{equation}
\end{theorem}

The theorem can also be stated requiring $u,v$, and $g$ to be defined in $\R$ and to be $2\pi$-periodic, or alternatively as a result on the circle.  Note that, as a consequence of \Cref{theorem:gy:non expansivity under J},  if $J(t)=|t|$ and the supremum of one of the functions $u$ and $v$ is less than or equal to the infimum of the other one, then equality in \eqref{rear:Euv} holds, whereas if the infima and suprema of both functions coincide, then \eqref{rear:|t|.sets} yields that all the superlevel sets of $u$ and $v$ are of the form $\bigcup_{k\in\Z}(I+2k\pi)$ for some  interval $I\in\R$ depending on the level.

It must be mentioned that most of the statements of \Cref{theorem:gy:non expansivity under J} are special cases of Burchard and Hajaiej \cite[Theorem 2]{Burchard Hajaiej}. However, the same comparing comments apply as those already mentioned right after \Cref{theorem:gy:J in R n-1}. In addition, we have the improvement that $J(0)=0$ is not assumed in \Cref{theorem:gy:non expansivity under J} (as it was in \Cref{theorem:gy:J in R n-1}), in contrary to \cite{Burchard Hajaiej} who do assume that $J(0)=0$. Removing this assumption, in particular in the case that $\inf J$ is not attained, is not straightforward.

As for the nonexpansivity \Cref{theorem:gy:J in R n-1}
 in $\R^n,$  \Cref{theorem:gy:non expansivity under J} is also a generalization of a Riesz rearrangement inequality, but now of a version on the circle. This Riesz rearrangement inequality is due to Baernstein \cite{Baernstein Circle} and is not as well known as the classical one in $\R^n$. 
There exists a generalization to the $n$-dimensional sphere, where in \eqref{eq:gy:intro:classical Riesz rearrangement} the space $\R^n$ is replaced by $\mathbb{S}^{n-1},$ the rearrangement is by geodesic balls around a chosen pole, and $g(x-y)$ is replaced by $K(x\cdot y),$ where $K:[-1,1]\to [0,+\infty)$ is a nondecreasing function.
In this form (i.e., with $g=g^{*\eper}$) the Riesz rearrangement inequality on the circle $\Sph^1$ was first proven independently by Baernstein and Taylor \cite{Baernstein Taylor} and by Friedberg and Luttinger \cite{FriedLutt}. They also generalized it to the sphere $\Sph^{n-1}$ with $n\geq2$ in \cite{Baernstein Taylor}, and on $\Sph^1$ but to a product of an arbitrary finite number of functions in \cite{FriedLutt}. We point out that these results cannot be deduced from the classical Riesz rearrangement inequality \eqref{eq:gy:intro:classical Riesz rearrangement} of the euclidean case.  It is still an open problem when $n\geq 2$ how a more general form of the Riesz rearrangement inequality on the sphere would look like ---i.e., not assuming $K$ to be symmetric and also rearranging it in some way; see \cite[Notes and Comments 8.11]{Baernstein book}. 
 A study of the cases of equality for the Riesz rearrangement inequality on the sphere can be found in  Burchard and Hajaiej \cite{Burchard Hajaiej}, respectively Baernstein \cite{Baernstein book}.
 We will gather a summary of these results on the circle in \Cref{thm:gy:sharp Friedberg Luttinger}.
 
Our proofs of the two nonexpansivity \Cref{theorem:gy:J in R n-1,theorem:gy:non expansivity under J} are very much inspired by \cite[Lemma A2]{Frank Seiringer}, which we generalize in the case of the Schwarz rearrangement to two functions $u$ and $v.$  Moreover, we give a slight simplification, since we do not need to consider the second derivative of $J$ in our proof, as is done in \cite{Frank Seiringer}. In short, the proof consists of  generalizing the known facts about the case $J(t)=t^2,$ i.e., the Riesz rearrangement inequality, to any convex function $J.$ We will prove \Cref{theorem:gy:non expansivity under J} in full detail in \Cref{section nonexpansive periodic}, since it is slightly harder and less known among the nonexpansivity results. Then, we will outline the proof of  \Cref{theorem:gy:J in R n-1} in \Cref{section nonexpansive cyl}, giving the details only where the proof differs from the periodic case.

\section{Some preliminaries on rearrangements}
\label{section preliminaries}

Let us start by recalling some basic facts on symmetrization and settle some notation. All sets and functions appearing in this section are assumed to be measurable with respect to the $n$-dimensional Lebesgue measure, where $n$ will vary by the context. We will frequently use the layer cake representation of a nonnegative function $u,$
$$
  u(x)=\int_0^{+\infty}dt\,\chi_{\{u>t\}}(x)=\sup\{t:\,x\in\{u>t\}\}.
$$

In this section we will deal with different types of rearrangements. We recall the definitions and most important properties that all these rearrangements have in common. To every measurable set $\Omega\subset\R^n$, $n\geq1$, we associate $\Omega^{\star},$ a set of same 
$n$-dimensional Lebesgue measure as $\Omega$, which we shall denote by $|\Omega|$, and with the following property: if $\Omega_1\subset\Omega_2$ then $\Omega_1^{\star}\subset\Omega_2^{\star}.$ For a characteristic function one defines $(\chi_{\Omega})^{\star}:=\chi_{\Omega^{\star}}.$  More generally, for every measurable function  $u$ defined in 
$\R^n$ or in a subset of it, we associate a rearranged function $u^{\star}$ by rearranging the superlevel sets of $|u|$, that is,
\begin{equation}
 \label{eq:gy:layer cake for general rearrange}
  u^{\star}(x):=\int_0^{+\infty}dt\,\chi_{\{|u|>t\}^{\star}}(x)=\sup\{t:\,x\in\{|u|>t\}^{\star}\}.
\end{equation}
This is well defined if $|\{|u|>t\}|<+\infty$ for all $t> 0.$
Hence $u^{\star}$ is nonnegative and the definition \eqref{eq:gy:layer cake for general rearrange} yields that, for every $\epsilon>0$, 
$$
  \{|u|>t+\epsilon\}^{\star}\subset\{u^{\star}>t\}\subset
  \{|u|>t\}^{\star}
  \text{ and, thus, }
  |\{|u|>t+\epsilon\}|\leq|\{u^{\star}>t\}|\leq
  |\{|u|>t\}|.
$$
This in turn easily gives that $u^{\star}$ is equimeasurable with $|u|$, i.e., $|\{u^{\star}>t\}|=|\{|u|>t\}|$ for all $t>0.$ Hence, up to a set of $n$-dimensional Lebesgue measure zero, it holds that
\begin{equation}\label{rear:lev.set.prop}
  \{u^{\star}>t\}=\{|u|>t\}^{\star}.
\end{equation}
Moreover, if for some measurable subset $\Omega\subset\R^n$ with $|\Omega|<+\infty$ it holds that 
\begin{equation}\label{eq:gy:u star cap Omega star}
  |\{u^{\star}>t\}\cap\Omega^{\star}|=|\{|u|>t\}\cap\Omega|\quad\text{ for all }t>0,
\end{equation}
and $F:\R\to\R$ is an absolutely continuous function, then
\begin{equation}
 \label{eq:gy:F abs contin ustar}
  \int_{\Omega}F(|u|)=\int_{\Omega^{\star}}F(u^{\star})
\end{equation}
---actually, this identity holds true for slightly more general functions $F$.
This follows immediately from Fubini's theorem, writing 
\begin{equation}
\begin{split}
  \int_{\Omega}F(|u|)&=\int_{\Omega}dx\,\Big(F(0)+\int_0^{|u(x)|}dt\,F'(t)\Big)\\
  &=F(0)|\Omega|+\int_{\Omega}dx\int_0^{+\infty}dt\,\chi_{(0,|u(x)|)}(t)F'(t)\\
  &= F(0)|\Omega|+\int_0^{+\infty}dt\,F'(t)|\{|u|>t\}\cap\Omega|.
\end{split}
\end{equation}
This argument also shows that \eqref{eq:gy:F abs contin ustar} holds in the case that $|\Omega|=+\infty$ if one further assumes, for example, that $F(0)=0$.

In the sequel $^\star$ will be one of the four rearrangements (or alternatively called symmetrizations) we will use in the different sections of this work, where $A\subset\R^n$ is a measurable set:
\begin{equation}
\begin{split}
   A^*:&\text{ Steiner symmetrization  in $\R^n$ with respect to the hyperplane $\{x_1=0\}.$}
  \\ A^{*\,n}:&\text{ Schwarz symmetrization in $\R^n$.}
  \\ A^{*\eper}:&\text{ periodic symmetrization in $\R^n$ with respect to the variable $x_1$.}
 \\  A^{{*\,n,1}}:&\text{ cylindrical symmetrization in $\R^n$ with respect to the variable $x_1$.}
\end{split}
\end{equation}
Clearly $A^*=A^{*\,1}.$ Be aware that we do not follow the notation of $(k,n)$-Steiner symmetrizations, where $1\leq k\leq n$ and $(1,n),$ $(n-1,n),$ and $(n,n)$ correspond to Steiner, cylindrical, and Schwarz symmetrization, respectively; see for instance \cite[Section 4]{Brock Solynin} where that notation was used.

We start with $A^*$, respectively its periodic version.
For a measurable set $A\subset\R$ with finite measure, its rearrangement is defined as the open interval $A^*:=\frac{1}{2}(-|A|,|A|)$, and the rearrangement of a function $u:\Omega\to\R,$ defined in some subset $\Omega\subset\R,$ is $u^*:\Omega^*\to\R$ defined by 
the layer cake formula \eqref{eq:gy:layer cake for general rearrange}.

For $A\subset\R^n$ we denote its section $A_{x'}:=\{x_1\in\R:\, (x_1,x')\in A\}\subset\R$ for every $x'\in\R^{n-1}.$ If $A\subset\R^n,$ then its Steiner symmetrization with respect to the hyperplane $\{x_1=0\}$ is defined by
$$
  A^{\ast}:=\bigcup_{x'\in\R^{n-1}}(A_{x'})^*\times\{x'\},
$$
with the agreement that if $|(A_{x'})^*|=0$ then 
$(A_{x'})^*\times\{x'\}=\emptyset$.
For a set $B\subset\R^n$ which is $2\pi$-periodic in the variable $x_1$, we define $B^{*\eper}$ by
$$
  B^{*\eper}:=\bigcup_{k\in\Z}\Big(\left(B\cap((-\pi,\pi)\times\R^{n-1})\right)^*+2k\pi e_1\Big),\quad e_1=(1,0,\ldots,0)\in\R^n.
$$
It follows from the definitions that, for every $A\subset\R^n$ and every $2\pi$-periodic set $B\subset\R^n$, 
\begin{equation}
 \label{eq:section and rearrange}
  (A_{x'})^*=(A^*)_{x'}\quad\text{and}\quad 
  (B_{x'})^{*\eper}=(B^{*\eper})_{x'}\quad \text{for all } x'\in\R^{n-1}.
\end{equation}

The periodic rearrangement of a measurable function $u:\R^{n}\to\R$ which is $2\pi$-periodic in the variable $x_1$ is defined by
$$
  u^{*\eper}(x):=\int_0^{+\infty}dt\,\chi_{\{|u|>t\}^{*\eper}}(x)
  =
  \sup\{t:\,x\in \{|u|>t\}^{*\eper}\}.
$$
It easily follows from \eqref{eq:section and rearrange} and the identity $\{|u(\cdot,x')|>t\}=\{|u|>t\}_{x'},$ that
\begin{equation}
 \label{eq:astper reduce to 1 dim}
  (u(\cdot,x'))^{*\eper}(x_1)=u^{*\eper}(x_1,x')\quad
  \text{for all $x_1\in\R$ and all $x'\in\R^{n-1}$},
\end{equation}
where the left hand side is simply the periodic rearrangement in $\R$ for every frozen value $x'.$

We now recall also the definition of the Schwarz symmetrization of a set. If $\Omega\subset\R^n$ is a measurable set of finite measure, then 
$$
  \Omega^{*\,n}:=\text{the open ball of $\R^n$ centered at the origin and of same measure as $\Omega.$}
$$

We will use the following property of rearrangements in a crucial way and, therefore, we state it in a lemma and provide its proof for the sake of completeness.

\begin{lemma}
\label{lemma:gy:increasing and low semicont}
Let $n\geq 1,$ $\Omega\subset\R^n$ be measurable, $u:\Omega\to\R$ be a measurable function, and $G:[0,+\infty)\to[0,+\infty)$ be a nondecreasing and lower semicontinuous function. Then,
\begin{equation}
 \label{eq:gy:lemma:G circ u ast}
  (G\circ |u|)^{*\,n}=G\circ (u^{\ast\,n})\quad\text{and}\quad
  (G\circ |u|)^{\ast\eper}=G\circ (u^{\ast\eper}),
\end{equation}
where for the second identity we assume that $n=1,$ $\Omega=\R,$ and $u$ is $2\pi$-periodic.
\end{lemma}

This lemma is sometimes stated with and sometimes without the assumption ``lower semincontinuous'' in the literature. Let us clarify this point.
If $G$ is nondecreasing, but not necessarily lower semicontinuous, then the identities \eqref{eq:gy:lemma:G circ u ast} still hold true for a.e.\! $x\in\Omega^*$, but not necessarily everywhere. Note that the rearrangement of a function is always lower semicontinuous because rearranged sets are \textit{open} balls. Consider, for instance, the following example: $\Omega=(-1,1),$ $u(x)=u^*(x)=1-|x|,$ and $G(t)=1$ if $|t|\geq 1/2$ and $0$ elsewhere. Then $(G\circ |u|)^*$ is the characteristic function of $(-1/2,1/2),$ but $G\circ (u^*)$ is the characteristic function of the closed interval $[-1/2,1/2].$ 

\begin{proof}[Proof of \Cref{lemma:gy:increasing and low semicont}]
In view of the layer cake representations, we have that
\begin{equation}
\begin{split}
  (G\circ |u|)^{*\,n}(x)&=\inf G+\int_{\inf G}^{\sup G} dt\,
  \chi_{\{G\circ |u|>t\}^{*\,n}}(x),\\
  G\circ (u^{*\,n})(x)&=\inf G+\int_{\inf G}^{\sup G} dt\,
  \chi_{\{G\circ (u^{*\,n})>t\}}(x).
\end{split}
\end{equation}
Therefore,
 it is enough to show that 
 $\{G\circ |u|>t\}^{*\,n}=\{G\circ (u^{*\,n})>t\}$ for all 
 $\inf G<t<\sup G$. Since $G$ is nondecreasing, lower semincontinuous, and $\inf G<t<\sup G$, the set $G^{-1}(t,+\infty)$ is equal to an open interval $(t_G,+\infty)$ for some $t_G>0$. Hence, it is enough to show that $\{|u|>t_G\}^{*\,n}=\{u^{*\,n}>t_G\},$ which is precisely \eqref{rear:lev.set.prop}. Note that, for $\star=*\,n$, the identity \eqref{rear:lev.set.prop} holds as an equality of sets (not only in measure) because both sides of the identity are open balls.
 This proves the first identity in \eqref{eq:gy:lemma:G circ u ast}.

The second identity in \eqref{eq:gy:lemma:G circ u ast}, which refers to the periodic rearrangement, follows from the first one applied to 
$\Omega=(-\pi,\pi)$. More precisely, since both functions 
$(G\circ |u|)^{\ast\eper}$ and $G\circ (u^{\ast\eper})$ are $2\pi$-periodic, it is enough to prove the identity in $(-\pi,\pi)$, and the latter follows from the fact that $f^*=f^{*\eper}$ in $(-\pi,\pi)$ for every $f.$
\end{proof}

We will also use the following rather elementary lemma to deal with sets of measure zero.

\begin{lemma}
\label{lemma:gy:slicing of layer cake}
Let $m\in\N$ and $f,\,g:\R^n\to[0,+\infty)$ be two nonnegative measurable functions.
\begin{itemize}
\item[$(i)$] Suppose that there exist a countable dense set 
$\{\tau_k\}_{k\in \N}\subset(0,+\infty)$ and a sequence of sets of zero $m$-dimensional Lebesgue measure $\{N_{k}\}_{k\in\N}$ such that, for every $k\in\N$, 
\begin{equation}
  \chi_{\{f>\tau_k\}}(y)=\chi_{\{g>\tau_k\}}(y)\quad\text{ for all }y \in
  \R^n\setminus N_{k}.
\end{equation}
Then $f=g$ almost everywhere in $\R^n$. 
In particular, if 
$\{f>\tau\}=\{g>\tau\}$ in measure for a.e.\! $\tau>0$, then $f=g$ almost everywhere in $\R^n$.

\item[$(ii)$] Suppose that $\{f>\tau\}$ is a ball (open or closed or containing some subsets of its boundary) for \textit{almost every} $\tau>0.$  Then, $\{f>\tau\}$ is a ball for \text{every} $\tau>0$ such that $|\{f>\tau\}|<+\infty$, in the sense that $\{f>\tau\}$ contains an open ball of the same measure.
\end{itemize}
\end{lemma}

Let us stress that, in \Cref{lemma:gy:slicing of layer cake}~$(i)$, we are assuming that the sets $\{f>\tau\}$ and $\{g>\tau\}$ only differ on a set of measure zero, but this set may depend on $\tau$. In this regard, if one actually knows that for a.e.\! $\tau>0$ the sets
$\{f>\tau\}$ and $\{g>\tau\}$ only differ on a set of measure zero independent of $\tau$, the fact that $f=g$ almost everywhere in $\R^n$ would follow directly from the layer cake decomposition of $f$ and $g$. However, the dependence on $\tau$ of the zero measure set makes the conclusion of \Cref{lemma:gy:slicing of layer cake}~$(i)$ not completely trivial.

Observe also that \Cref{lemma:gy:slicing of layer cake}~$(ii)$ may not hold if $|\{f>\tau\}|<+\infty$ is not assumed. For example, think on the balls $B_n(ne_1)$, whose boundary passes through the origin for every $n.$ Then, the union of the balls over all $n$ gives a half space.

\begin{proof}[Proof of \Cref{lemma:gy:slicing of layer cake}] 
Let us first address the proof of $(i)$. We will see that $|\{|f-g|>0\}|=0$. Note that
\begin{equation}
\begin{split}
\{|f-g|>0\}=\bigcup_{j=1}^{+\infty}\{|f-g|>1/j\}
=\bigcup_{j=1}^{+\infty}\Big(\{f-g>1/j\}\cup\{g-f>1/j\}\Big).
\end{split}
\end{equation} Hence, it is enough to check that 
$|\{f-g>1/j\}|=|\{g-f>1/j\}|=0$ for all $j\geq1$. We will prove that $|\{f-g>1/j\}|=0$; the case of $|\{g-f>1/j\}|$ follows analogously.

Assume that $x\in \{f-g>1/j\}$, which means that $f(x)-g(x)>1/j$. Since $\{\tau_k\}_{k\in \N}$ is dense in $(0,+\infty)$, there exists $k_0\in\N$ such that $g(x)\leq\tau_{k_0}<f(x)$. Therefore,
$x\in\{f>\tau_{k_0}\}\cap\{g\leq\tau_{k_0}\}
=\{f>\tau_{k_0}\}\setminus\{g>\tau_{k_0}\}$. From this we deduce that
$$\{f-g>1/j\}\subset\bigcup_{k\in\N}\Big(\{f>\tau_{k}\}\setminus\{g>\tau_{k}\}\Big).$$
Now, the fact that $\chi_{\{f>\tau_k\}}(y)=\chi_{\{g>\tau_k\}}(y)$ for all $y \in \R^n\setminus N_{k}$ and  that $|N_k|=0$ yields  
$|\{f>\tau_{k}\}\setminus\{g>\tau_{k}\}|=0$
for all $k\in\N$. Thus, $|\{f-g>1/j\}|=0$ for all $j\geq1$, as desired.

We now prove $(ii)$. Let $\tau_i\downarrow \tau$ be a monotone sequence such that $\{f>\tau_i\}$ are balls with radii $R_i$. From the inclusions $\{f>\tau_i\}\subset \{f>\tau_{i+1}\},$ we get that the sequence $\{R_i\}_i$ is monotone nondecreasing and has a limit $R_{\infty}$, which  is finite because $|\{f>\tau\}|<+\infty$ by assumption. The centers of the balls $\{f>\tau_i\}$ must have a converging subsequence, and it can be easily seen that the limit of such a sequence must be unique. One then concludes the proof of $(ii)$ simply using that $\{f>\tau\}=\bigcup_{i=1}^{+\infty}\{f>\tau_i\}.$ 
\end{proof}

\section{Nonexpansivity under periodic rearrangement}
\label{section nonexpansive periodic}

The aim of this section is to prove  \Cref{theorem:gy:non expansivity under J}.
The proof is based on the following version of the Riesz rearrangement inequality on the circle.


\begin{theorem}
[\cite{Baernstein Circle,Baernstein book,Baernstein Taylor,Burchard Hajaiej,FriedLutt}]\label{thm:gy:sharp Friedberg Luttinger}
Let $f,h:(-\pi,\pi)\to \R$ be two nonnegative measurable functions and $g:(-2\pi,2\pi)\to [0,+\infty)$ be a $2\pi$-periodic measurable function. 
Then, 
\begin{equation}\label{eq:thm FriedLutt inequality}
    \int_{-\pi}^{\pi}dx\int_{-\pi}^{\pi} dy\,f(x)g(x-y)h(y)
    \leq 
    \int_{-\pi}^{\pi}dx\int_{-\pi}^{\pi}dy\, f^{*}(x)g^{*\eper}(x-y)h^{*}(y).
\end{equation}
If in addition $g=g^{*\eper}$, $g$ is  decreasing in $(0,\pi),$ 
and the left-hand side of \eqref{eq:thm FriedLutt inequality} is finite,
then
equality holds in \eqref{eq:thm FriedLutt inequality} if and only if at least one of the following conditions holds:  

\begin{itemize}
\item[$(i)$] One of the two functions $f$ or $h$ is constant almost everywhere (and the other one can be anything).

\item[$(ii)$] If $f$ and $h$ are extended to $\R$ in a $2\pi$-periodic way, then there exists $z\in \R$ such that $f(x)=f^{*\eper}(x+z)$ and $h(x)=h^{*\eper}(x+z)$ for almost every $x\in\R$.
\end{itemize}
\end{theorem}

The inequality \eqref{eq:thm FriedLutt inequality} in the case $g=g^{*\eper}$ was first discovered, independently, in \cite{Baernstein Taylor} and \cite{FriedLutt}. Both references contain more general inequalities: \cite{Baernstein Taylor} deals with the sphere $\mathbb{S}^n$, whereas \cite{FriedLutt} deals with a product of more than three functions defined in~$\R$. 
The result \cite[Theorem 7.3]{Baernstein book} is also more general than  \Cref{thm:gy:sharp Friedberg Luttinger}, as it deals with a Riesz rearrangement inequality on the sphere $\mathbb{S}^n$ and not only on the circle. Moreover, \cite{Baernstein book} treats more general functions of $f(x)$ and $h(y)$ than simply the product $f(x)h(y).$ 
Without the assumption $g=g^{*\eper}$, the inequality \eqref{eq:thm FriedLutt inequality} can be found in \cite{Baernstein Circle}. We mention that we stated the inequality in this general form, but we will actually only use the case when $g=g^{*\eper}$.

Instead, the statement in \Cref{thm:gy:sharp Friedberg Luttinger} concerning equality in \eqref{eq:thm FriedLutt inequality} follows from Burchard and Hajaiej \cite[Theorem 2]{Burchard Hajaiej}, who treated the case of equality in $\mathbb{S}^n$ for the first time. For this result, we also cite \cite[Theorem 7.3]{Baernstein book} since, being less general than \cite{Burchard Hajaiej}, fits precisely with our setting.

Let us finally mention that we find \cite{Baernstein book} to be the simplest reference for looking up all the statements of \Cref{thm:gy:sharp Friedberg Luttinger} that we need in the present work.

For the proof of  \Cref{theorem:gy:non expansivity under J}, it will be convenient to introduce the following abbreviation
\begin{equation}\label{energy.E.J}
  E[u,v,g]:=\int_{-\pi}^{\pi}dx\int_{-\pi}^{\pi}dy\,J(u(x)-v(y))g(x-y),
\end{equation}
where $J:\R\to [0,+\infty)$ is some convex function.
Moreover, whenever we will assume that $g=g^{*\eper}$, as is the setting when studying cases of equality, as well as in all  our applications, we will omit that argument and will simply use the notation
$$
  E[u,v]:=E[u,v,g]\quad\text{ if }g=g^{*\eper}.
$$
 The main idea of the proof is based on that of the nonexpansivity of the $L^p$ norm under rearrangement without the convolution kernel $g$ (see \cite[Theorem 3.5]{Lieb Loss}), but one uses in the proof the periodic Riesz rearrangement inequality of  \Cref{thm:gy:sharp Friedberg Luttinger} instead of the Hardy-Littlewood inequality $\int fh\leq \int f^* h^*.$ This same strategy of the proof was also used by Frank and Seiringer \cite[Appendix A, Lemma A.2]{Frank Seiringer}.

Before starting the proof let us make some comments on some of the hypothesis' of the theorem.
Note that, contrary to the case $u=v,$ strict convexity of $J$ in $[0,+\infty)$ or $(-\infty,0]$ is not enough to characterize cases of equality, as in \cite[first paragraph in the proof of Lemma~A.2]{Frank Seiringer}. Consider, for instance, a function $J$ which is identically zero in the negative real line and strictly convex in $[0,+\infty).$ Then, for any pair of nonnegative functions $u$ and $v$ such that $\inf v\geq \sup u$ it holds that $E[u,v]=0=E[u^*,v^*].$

\begin{remark}
\label{remark:negative functions}
{\em
If $J$ is an even, nonnegative, and convex function with its minimum at $0$, then the inequality in \Cref{theorem:gy:non expansivity under J} remains true even if we drop the assumption on the sign of $u$ and $v$, since $|u-v|\geq ||u|-|v||$, but one would have to ensure by some other assumption than $|\{u>\tau\}|<+\infty$ that the rearrangement is well defined. For instance, if $J(t)=t^2$ then the inequality still holds, as can be easily seen by the factorization $(u-v)^2=u^2+v^2-2uv$, the fact that $u\leq |u|$, and \Cref{thm:gy:sharp Friedberg Luttinger}. That $(ii)$ implies the equality in \eqref{rear:Euv} obviously remains true in this special case, but that $(i)$ implies the equality in \eqref{rear:Euv} does not hold if we do not assume $u$ and $v$ to be nonnegative, even for $J(t)=t^2$. Consider, for instance, that $u= c>0$ is some positive constant, in which case 
\begin{align}
  E[u,v]=&\int_{-\pi}^{\pi}dx\,g(x)\int_{-\pi}^{\pi}dy\,(c-v(y))^2
  =\int_{-\pi}^{\pi}dx\,g(x)\int_{-\pi}^{\pi}dy\,(c^2-2cv(y)+v^2(y))
  \\
  =& \int_{-\pi}^{\pi}dx\,g(x)\int_{-\pi}^{\pi}dy\,(c^2-2cv(y)+(v^*)^2(y)).
\end{align}
Now, for any $v\geq 0$ not necessarily symmetric, we have 
$\int v=\int |v|=\int v^*$. However, if $v$ takes on negative values on a set of positive measure, then $\int v<\int |v|=\int v^*,$ and therefore $E[u^*,v^*]$ is strictly smaller than $ E[u,v]$. Thus, the condition of $(i)$ is satisfied, since $u$ is constant, but there is no equality in \eqref{rear:Euv}.}
\end{remark}

\begin{proof}[Proof of  \Cref{theorem:gy:non expansivity under J}.]
The proof of the theorem will be divided in five steps. We will first assume that $\min J$ exists to carry out the proof of the theorem for this case in Steps~1 to~4. The proof of the case that 
$\inf J$ is not attained is slightly different and is given in Step 5. Recall that we write $E[u,v]$ instead of $E[u,v,g]$ whenever $g=g^{*\eper}.$
\medskip

\noindent\underline{{\em Step 1}  (preliminaries):} 
\smallskip

The purpose of this first step is to show that we can assume that $J(0)=0$. To see this, let $J(t_0)=\min J$ and consider the new  function $\widetilde{J}(t)=J(t+t_0)-\min J,$ which vanishes at $t=0$, is nonnegative, and convex. Denote the corresponding functional for $\widetilde{J}$ by $\widetilde{E}[u,v,g].$ The term in $\widetilde{E}[u,v,g]$ involving $\min J$ is equal to
\begin{equation}\label{eq:minJ_unnecessary}
  -\min J\int_{-\pi}^{\pi}dx\int_{-\pi}^{\pi}dy\,g(x-y)
  =-2\pi \min J\, \|g\|_{L^1(-\pi,\pi)}
  =
  -2\pi \min J\, \|g^{*\eper}\|_{L^1(-\pi,\pi)},
\end{equation}
which does not depend on $u$ and $v$. Therefore, this term will have no effect when showing the inequality or studying the case of equality.

Let us first show that it is sufficient to prove the inequality for $\widetilde{E}$ to get it also for $E$. Assume that the inequality holds for $\widetilde{J}$ and any two nonnegative functions $u$ and $v$. On the one hand, if $t_0\geq0,$ by \eqref{eq:minJ_unnecessary} we have
\begin{equation}
\begin{split}
  E[u,v,g]
  &=\intp dx\intp dy\,\widetilde{J}(u(x)-(v(y)+t_0))g(x-y)
  +2\pi \min J\, \|g\|_{L^1(-\pi,\pi)}\\
  &=\widetilde{E}[u,v+t_0,g] +2\pi \min J\, \|g\|_{L^1(-\pi,\pi)}\\
  &\geq  \widetilde E[u^*,(v+t_0)^*,g^{*\eper}]
  +2\pi \min J\, \|g^{*\eper}\|_{L^1(-\pi,\pi)}\\
  &=  \widetilde E[u^*,v^*+t_0,g^{*\eper}]
  +2\pi \min J\, \|g^{*\eper}\|_{L^1(-\pi,\pi)}\\
  &=\intp dx\intp dy\, \widetilde{J}(u^*(x)-(v^*(y)+t_0))g^{*\eper}(x-y)+2\pi \min J\, \|g^{*\eper}\|_{L^1(-\pi,\pi)}\\
  &=E[u^*,v^*,g^{*\eper}],
\end{split}
\end{equation}
where we have used that $v,t_0\geq 0$ and, hence, $v+t_0$ is nonnegative and $(v+t_0)^*=v^*+t_0.$  On the other hand, if $t_0<0$ one proceeds in a similar way, but using now the inequality for $\widetilde{J}$ and the nonnegative functions $u-t_0$ and $v$.

Let us now assume that $g=g^{*\eper},$ that $E[u,v]=E[u^*,v^*],$ and that we can characterize the cases of equality for $\widetilde{E}.$ If $t_0\geq0$ then, using that
\begin{equation}
\begin{split}
  \widetilde{E}[u,v+t_0]
  &=E[u,v]-2\pi \min J\, \|g\|_{L^1(-\pi,\pi)}\\
  &=E[u^*,v^*]-2\pi \min J\, \|g^{*\eper}\|_{L^1(-\pi,\pi)}
  =\widetilde{E}[u^*,(v+t_0)^*]
\end{split}
\end{equation}
and that $(v+t_0)^*=v^*+t_0$, the cases of equality for $E$ are also characterized. As before, if $t_0<0$ simply consider $u-t_0$, for which $(u-t_0)^*=u^*-t_0$. 

All these considerations show that it is enough to prove the theorem for $\widetilde J$ or, equivalently, that we can assume that $J(0)=0$ without loss of generality.

\medskip

\noindent\underline{{\em Step 2} (proof of the inequality):} 
\smallskip 

In this section we will prove \eqref{rear:Euv}. In view of Step 1, to do it we can assume that $\min J=J(0)=0$.
Let us decompose $J=J_{+}+J_{-},$ where $J_+(t):=J(t)$ for $t\geq 0$ and $J_+(t):=0$ for $t<0$, and define the corresponding functionals
\begin{equation}\label{def.E_pm}
  E_{\pm}[u,v,g]:=\intp dx\intp dy\, J_{\pm}(u(x)-v(y))g(x-y).
\end{equation}
Note that, by the assumptions on $J$, $J_\pm$ are nonnegative convex functions with $J_\pm(0)=0$. In particular, $J_+$ is nondecreasing in $\R$.

Using this decomposition of $J$, it is clear that it is enough to prove the inequality for $E_{\pm}$ to get it also for $E$. Actually, we only need to show that $E_+$ is nonincreasing under rearrangement. Once this is proven, it easily follows that $E_-$ is nonincreasing under rearrangement too by applying the proof for done for $E_+$ to $t\mapsto J_-(-t)$ (which is also nonnegative, convex, and vanishes in $(-\infty,0]$)
and interchanging the roles of $u$ and $v.$ Therefore, from now on we will only focus on the proof of \eqref{rear:Euv} for $J_+$.

Note that $J_+'$ is nondecreasing, nonnegative, and continuous almost everywhere (see, for instance, \cite[Theorem 25.3]{Rockafellar}) and hence, redefining $J'$ on  a set of measure zero, we can assume that $J_+'$ is lower semicontinuous in $\R$. We chose this representative of $J_+'$ because, then, any translation of $J_+'$ satisfies the assumptions of the function $G$ in \Cref{lemma:gy:increasing and low semicont}. Also, recalling that $J_+=0$ in $(-\infty,0]$, we have $J_+'(0)=0.$

Since $J_+$ is locally Lipschitz continuous, it is also absolutely continuous on any bounded set, and 
therefore
\begin{equation}
  J_+(u(x)-v(y))
  =
  \int_{u(x)}^{v(y)}d\tau\,\frac{d}{d\tau}\Big(J_+(u(x)-\tau)\Big)
  =
  -\int_{u(x)}^{v(y)}d\tau\,J_+'(u(x)-\tau).
\end{equation}
Recall that $J'_+= 0$ in $(-\infty,0],$ thus
\begin{equation}
\begin{split}
  J_+(u(x)-v(y))
  &=
  \int_{v(y)}^{+\infty}d\tau\,J_+'(u(x)-\tau)
  =
  \int_0^{+\infty}d\tau\,J_+'(u(x)-\tau)\chi_{\{v{\leq}\tau\}}(y)
  \\
  &=
  \int_0^{+\infty}d\tau\,J_+'(u(x)-\tau)(1-\chi_{\{v>\tau\}}(y)).
\end{split}
\end{equation}
Since $E_+[u,v,g]$ is finite, an application of Fubini's theorem gives
\begin{equation}
 \label{eq:gy:E u v finite and Fubini}
  E_+[u,v,g]=\int_0^{+\infty}d\tau\left(\int_{-\pi}^{\pi}dx\int_{-\pi}^{\pi}dy\,
  J_+'(u(x)-\tau)g(x-y)(1-\chi_{\{v>\tau\}}(y))\right).
\end{equation}

We will now distinguish first a special case and then we will treat the general case.

\medskip
\noindent\underline{{\em Step 2.1} (proof of the inequality for $u$ or $|J'|$ bounded):} 
\smallskip 

Note first that if $u$ is bounded, that is, if $\esssup u\leq M<+\infty$ for some $M>0,$ then we can estimate, for any $\tau\geq 0$,
\begin{equation}
\begin{split}
  \intp dx\intp dy\, J_+'(u(x)-\tau)g(x-y)
  &\leq  2\pi \|g\|_{L^1(-\pi,\pi)}J_+'(M-\tau)<+\infty. 
\end{split}
\end{equation}
Similarly, if $|J'|$ is bounded (which occurs, for instance, if $J(t)=|t|$), the double integral on the left hand side is finite too.
Therefore, in any case we can split $E_+[u,v,g]$ given in \eqref{eq:gy:E u v finite and Fubini} as
\begin{equation}
 \label{eq:gy:A and B.E_+}  
 E_+[u,v,g]=\int_0^{+\infty}d\tau\,\big(A(u,g,\tau)-B(u,v,g,\tau)\big),
\end{equation}
where $A$ and $B$ are both finite and nonnegative for every 
$\tau,$ and are given by
\begin{equation}
 \label{eq:gy:A and B}
  \begin{split}
  A(u,g,\tau)&:=\int_{-\pi}^{\pi}dx\int_{-\pi}^{\pi}dy\,
  J_+'(u(x)-\tau)g(x-y),
   \\
  B(u,v,g,\tau)&:=\int_{-\pi}^{\pi}dx\int_{-\pi}^{\pi}dy\,
  J_+'(u(x)-\tau)g(x-y)\chi_{\{v>\tau\}}(y).
 \end{split}
\end{equation} 
Observe that $A(u,g,\tau)\geq B(u,v,g,\tau)$, hence the integrand on the right hand side of \eqref{eq:gy:A and B.E_+} is nonnegative. 

We will show  that, for every $\tau\geq 0$, the functional $A$ is unaltered under the rearrangement of $u,v$, and $g$, and that $B$ does not decrease. Once we have shown this, the inequality of the theorem is shown as follows, using the abbreviations $A(\tau)\equiv A(u,g,\tau)$ and $A^*(\tau)\equiv A(u^*,g^{*\eper},\tau)$, and analogously for $B$:
\begin{equation}\label{eq:gy:E u v split in A B and ineq}
\begin{split}
  E_+[u,v,g]&= \int_0^{+\infty}d\tau\,(A(\tau)-B(\tau))
  =\int_0^{+\infty}d\tau\,(A^*(\tau)-B(\tau))
  \\
  &=\int_0^{+\infty}d\tau\,(A^*(\tau)-B^*(\tau)+B^*(\tau)-B(\tau))
   \\
  &=E_+[u^*,v^*,g^{*\eper}]+\int_0^{+\infty}d\tau\, (B^*(\tau)-B(\tau))
  \geq E_+[u^*,v^*,g^{*\eper}].
 \end{split} 
\end{equation}
Observe that one cannot assume in general that $\int d\tau\, A(\tau)$ and $\int d\tau\, B(\tau)$ are finite.\footnote{This would be essentially equivalent to assuming, for instance if $J(t)=|t|$ and $g\equiv 1$, that $u\in L^1(-\pi,\pi)$, since $\int d\tau A(\tau)=2\pi \|u\|_{L^1}$.} In principle, only the difference $A(\tau)-B(\tau)$ is known to be integrable in $\tau\in(0,+\infty)$. 

It only remains to show the claimed behaviors of $A$ and $B$ under rearrangement.
We first deal with $A.$ Since $g$ is $2\pi$-periodic and nonnegative, we have $\int_{-\pi}^{\pi}dy\,g(x-y)=\|g\|_{L^1(-\pi,\pi)}
=\|g^{*\eper}\|_{L^1(-\pi,\pi)}$ for all $x$. Hence, it is sufficient to show that
\begin{equation}
 \label{eq:J prime u ast}
  \int_{-\pi}^{\pi}dx\,J_+'(u(x)-\tau)=\int_{-\pi}^{\pi}dx\,J_+'(u^*(x)-\tau).
\end{equation}
Using that $\intp f=\intp f^*$ for the nonnegative function $f:=J_+'\circ(u-\tau)$ we can first replace the integrand on the left hand side of \eqref{eq:J prime u ast} by $(J_+'\circ(u-\tau))^*$, and then apply  \Cref{lemma:gy:increasing and low semicont} with $G(s):=J_+'(s-\tau)$ for all $s\geq0$, since $u=|u|$ by assumption. We therefore get
\begin{equation}
 \label{eq:J and ast interchang}
  (J_+'(u-\tau))^*(x)
  =
  J_+'(u^*(x)-\tau),
\end{equation}
and  \eqref{eq:J prime u ast} follows.

The claim that $B$ does not decrease under rearrangement follows from  \Cref{thm:gy:sharp Friedberg Luttinger} using \eqref{eq:J and ast interchang} and \eqref{rear:lev.set.prop}.

\medskip
\noindent\underline{{\em Step 2.2} (proof of the inequality for general $u$):} 
\smallskip

We will now remove the assumption that $u$ is bounded. To do it, for $M>0$ define 
\begin{equation}
 \label{eq:gy:u M Def}
  u_M:=\min\{u,M\}\quad\text{and}\quad
  (u^*)_M:=\min\{u^*,M\}.
\end{equation}
Observe that if we replace in the integrand of  \eqref{def.E_pm}, in the case $J_+$, the function $u$ by $u_M$, then this integrand is nondecreasing for every $(x,y)$ as $M$ increases. Hence the monotone convergence theorem gives that
$$
  \lim_{M\uparrow+\infty}E_+[u_M,v,g]=E_+[u,v,g]\quad
  \text{ and }\quad
  \lim_{M\uparrow+\infty}E_+[(u^*)_M,v^*,g^{*\eper}]=E_+[u^*,v^*,g^{*\eper}].
$$
The proof of the inequality for general $u$ follows from Step 2.1 and the fact that $(u^*)_M=(u_M)^{*}.$

\medskip

\noindent\underline{{\em Step 3} (case of equality for $J$ strictly convex):} 
\smallskip 

Throughout this step we will assume that $J$ is strictly convex, that $g=g^{*\eper}$ is decreasing in $(0,\pi),$ and shall write $E[u,v]$ instead of $E[u,v,g].$ The purpose of this step is to show that there is equality in \eqref{rear:Euv} if and only if $(i)$ and/or $(ii)$ in \Cref{theorem:gy:non expansivity under J} holds.
 \medskip
 
\noindent\underline{{\em Step 3.1} (case of equality for $u$ or $v$ constant, or $u$ and $v$ as in \Cref{theorem:gy:non expansivity under J}~$(ii)$):} 
\smallskip

Note that $(ii)$ trivially yields  equality in \eqref{rear:Euv}. 

We will show here that if one of the functions $u$ or $v$ is constant (that is, if \Cref{theorem:gy:non expansivity under J}~$(i)$ holds), then there is equality in \eqref{rear:Euv}. Let us assume that $u(x)=c\geq 0$ for all $x$; the case of $v$ constant is analogous. Then, by the $2\pi$-periodicity of $g,$ 
$$
  E[u,v]=\int_{-\pi}^{\pi}dy\,J(c-v(y))\left(\int_{-\pi}^{\pi}dx\,g(x-y)\right)dy
  =
  \|g\|_{L^1(-\pi,\pi)}\int_{-\pi}^{\pi}dy\,J(c-v(y)).
$$
Now, by the equimeasurability indentity \eqref{eq:gy:F abs contin ustar} applied to $F(t):=J(c-t)$ we get that
\begin{align}
  \int_{-\pi}^{\pi}dy\,J(c-v(y))=\int_{-\pi}^{\pi}dy\,J(c-v^*(y)).
\end{align}
This shows that $E[u,v]=E[u^*,v^*].$

\medskip
\noindent\underline{{\em Step 3.2} (case of equality for $u$ and $v$ bounded):} 
\smallskip

In this step we will show that if both $u$ and $v$ are bounded, neither $u$ nor $v$ are constant, and there is equality in \eqref{rear:Euv}, then \Cref{theorem:gy:non expansivity under J}~$(ii)$ holds.
 
Set $S:=\operatorname{ess}\sup v$ and 
$$
  I_S(t):=\left\{\begin{array}{rl}
                J(t-S)-tJ'(-S)-J(-S) & \text{ if }t\geq 0,
                \\
                0 & \text{ if }t<0.
               \end{array}\right. 
$$
If $J$ is not differentiable at $-S,$ then take instead of $J'(-S)$ in the above formula any value of a subgradient (i.e., the slope of a tangent touching the graph of $J$ at $-S$). Note that $I_S(0)=0.$ Therefore, $I_S$ is a convex and nonnegative function, and strictly convex in $[0,+\infty).$ Also, observe that $u+S-v\geq0$, thus 
\begin{equation}\label{def.I_S.J}
J(u(x)-v(y))=I_S(u(x)+S-v(y))
+(u(x)+S-v(y))J'(-S)+J(-S)
\end{equation}
for all $x,y$. Regarding the term $(u(x)+S-v(y))J'(-S)$,
by the $2\pi$-periodicity of $g$ we have that (recall also that we are assuming $u$ and $v$ bounded)\footnote{
Note that this step will not work in $\R$ instead of $(-\pi,\pi),$ since the integrals may be infinite ---compare with the proof of \Cref{theorem:gy:J in R n-1} in \Cref{section nonexpansive cyl}; see in particular the paragraph after \eqref{eq:gy:u ast m minut tau}.}
\begin{equation}
\begin{split}
 \intp dx\intp dy\, (u(x)+S-v(y))g(x-y)
 =\|g\|_{L^1(-\pi,\pi)}\Big(\int_{-\pi}^{\pi}(u+S)
 -\int_{-\pi}^{\pi}v\Big).
\end{split}
\end{equation}
Hence, this expression coming from $tJ'(-S)$ with $t=u(x)+S-v(y)$
in the definition of $I_S(t)$ is unaltered by the rearrangement of $u$ and $v$ because $(u+S)^*=u^*+S$. The same conclusion holds for the expression coming from $J(-S)$ in the definition of $I_S(t)$. Therefore, 
defining $u_S:=u+S$, the assumption of equality in \eqref{rear:Euv} and \eqref{def.I_S.J} yield 
\begin{equation}\label{E_S.u.ustar}
\begin{split}
  E_S[u,v]&:=\intp dx\intp dy\,  I_S(u_S(x)-v(y))g(x-y)
   \\
  &=\intp dx\intp dy\, I_S((u_S)^*(x)-v^*(y))g(x-y)
  =E_S[u^*,v^*],
  \end{split}
\end{equation}
where we have used that $(u_S)^*=(u^*)_S$. Since $I_S(t)=0$ for $t\leq0,$  we can proceed exactly as we did in Step 2 to obtain \eqref{eq:gy:E u v finite and Fubini} and write
\begin{equation}\label{E_S.int.tau.}
\begin{split}
    E_S[u,v]=\int_0^{+\infty}d\tau\left(\int_{-\pi}^{\pi}dx\int_{-\pi}^{\pi}dy\,
  I_S'(u_S(x)-\tau)g(x-y)(1-\chi_{\{v>\tau\}}(y))\right).
\end{split}
\end{equation}
Since we assumed that $u$ is bounded, $u_S$ too. Thus, all conclusions of Step 2 hold with $E_+$ replaced by $E_S$ and $A,B$ replaced by $A_S,B_S$ where
\begin{equation}
  \begin{split}
  A_S(u_S,\tau)&:=\int_{-\pi}^{\pi}dx\int_{-\pi}^{\pi}dy\,
  I_S'(u_S(x)-\tau)g(x-y), \\
  B_S(u_S,v,\tau)&:=\int_{-\pi}^{\pi}dx\int_{-\pi}^{\pi}dy\,
  I_S'(u_S(x)-\tau)g(x-y)\chi_{\{v>\tau\}}(y).
 \end{split}
\end{equation} 
In particular, $A_S$ is unaltered under rearrangement and the analog of \eqref{eq:gy:E u v split in A B and ineq} for $I_S'$ and $u_S$ holds. Therefore, using that the integrand of $B_S$ is nonnegative, that $B_S$ does not decrease under rearrangement, \eqref{E_S.u.ustar}, and \eqref{E_S.int.tau.}, we must have that $B_S(u_S,v,\tau)=B_S((u_S)^*,v^*,\tau),$ i.e., 
\begin{equation}\label{eq:eqIS_uv_nonconst}
  \begin{split}
    \intp dx \intp dy\, I_S'(u_S&(x)-\tau)g(x-y)\chi_{\{v>\tau\}}(y)\\
  &=
  \intp dx \intp dy\, (I_S'(u_S-\tau))^*(x)g(x-y)\chi_{\{v^*>\tau\}}(y)
\end{split}
\end{equation}
for almost every $\tau>0$ (here we have also used \eqref{eq:J and ast interchang} with $I_S'$ instead of $J_+'$). We will therefore be in position to apply the conclusions of  \Cref{thm:gy:sharp Friedberg Luttinger}, as done below.

We now use the assumption that $u$ and $v$ are not constant almost everywhere. Firstly, this means that 
$$
  \chi_{\{v>\tau\}}(y)\text{ is not constant in $y$ for any }\tau\in (\operatorname{ess}\inf v,\operatorname{ess}\sup v)= (\operatorname{ess}\inf v,S)
  \neq\emptyset.
$$
For this range of $\tau$ we also obtain that 
\begin{equation}\label{eq:uS minus tau}
  u_S(x)-\tau\geq S-\tau>0.
\end{equation}Since $I_S$ is strictly convex in $[0,+\infty),$ $I_S'$ is  increasing and, therefore, one-to-one in the positive real line. Thus,
$$
  I_S'(u(x)-\tau)\text{ is not constant in $x$ for any }\tau\in (\operatorname{ess}\inf v,S).
$$
Assuming that $u$ and $v$ are extended to $\R$ in a $2\pi$-periodic way, from \eqref{eq:eqIS_uv_nonconst} and  the case of equality in  \Cref{thm:gy:sharp Friedberg Luttinger} we deduce that for almost every $\tau\in(\operatorname{ess}\inf v,S)$ there exists $z_{\tau}\in\R$ such that
\begin{equation}\label{eq:gy:chi v bigger tau a.e.y }
\begin{split}
   I_S'(u_S(x)-\tau)&=    \left(I_S'(u_S-\tau)\right)^{*\eper}(x+z_{\tau})=I_S'((u_S)^{*\eper}(x+z_{\tau})-\tau),\\
   \chi_{\{v>\tau\}}(y)&=\chi_{\{v^{*\eper}>\tau\}}(y+z_{\tau})
\end{split}
\end{equation}
for almost every $x,y\in\R$, 
where we have also used  \Cref{lemma:gy:increasing and low semicont} as in \eqref{eq:J and ast interchang} in the first line in \eqref{eq:gy:chi v bigger tau a.e.y }.
Since $I_S'$ is one-to-one in the positive real line, the first identity is equivalent to
\begin{equation}
 \label{eq:uS.uSper.tau}
  u_S(x)=(u_S)^{*\eper}(x+z_{\tau})\text{ for a.e.\! }x\in\R,\text{ for a.e.\! } \tau\in  (\operatorname{ess}\inf v,S).
\end{equation}
From this it follows that $z_{\tau}$ is indeed independent of $\tau$ because $u$ is not constant. In more detail, that $u$ is not constant yields that $(u_S)^{*\eper}$ is a nontrivial symmetric decreasing function in $(-\pi,\pi)$, which in particular yields that $(u_S)^{*\eper}(\cdot+z_{\tau_1})=(u_S)^{*\eper}(\cdot+z_{\tau_2})$ almost everywhere if and only if $z_{\tau_1}=z_{\tau_2}$ modulo $2\pi$, since $(u_S)^{*\eper}$ has minimal period $2\pi$. 

Finally, using that $(u_S)^{*\eper}=u^{*\eper}+S$, \eqref{eq:uS.uSper.tau}, and the fact that $z_\tau$ is indeed independent of $\tau$, we conclude that $u(x)=u^{*\eper}(x+z)$ for some $z\in\R$ and almost every $x\in\R$.
In addition, from \eqref{eq:gy:chi v bigger tau a.e.y } we now obtain that
\begin{equation}
 \label{eq:gy:chi v bigger tau y}
   \chi_{\{v>\tau\}}(y)=\chi_{\{v^{*\eper}>\tau\}}(y+z)
   \text{ for a.e.\! }y\in\R,\text{ for a.e.\! } \tau\in  (\operatorname{ess}\inf v,\operatorname{ess}\sup v).
\end{equation}
This clearly holds also for all $\tau\in(0,+\infty)\setminus(\operatorname{ess}\inf v,\operatorname{ess}\sup v)$ by definition of the rearrangement. Thus, using  \Cref{lemma:gy:slicing of layer cake}~$(i)$ we get that $v(y)=v^{*\eper}(y+z)$ for almost every $y\in\R.$

\medskip
\noindent\underline{{\em Step 3.3} (case of equality for general $u$ and $v$):} 
\smallskip

To complete the proof of Step 3, it only remains to remove the boundedness assumption in Step 3.2. This is precisely what we will do here. More precisely, we will show that if neither $u$ nor $v$ are constant and there is equality in \eqref{rear:Euv}, then \Cref{theorem:gy:non expansivity under J}~$(ii)$ holds. 

Observe first that the equality $E[u,v]=E[u^*,v^*]$ forces the validity of both equalities
\begin{equation}
 \label{eq:gy:E plus u v is and E minus}
  E_+[u,v]=E_+[u^*,v^*]\quad\text{ and }\quad E_-[u,v]=E_-[u^*,v^*],
\end{equation}
where $E_\pm$ have been introduced in \eqref{def.E_pm}.
This is because $E=E_+ + E_-$, $E_+$ and $E_-$ are nonnegative, and both are nonincreasing under rearrangement as shown in Step 2.

Let $u_M$ and $v_M$ be defined as in \eqref{eq:gy:u M Def}.
Observe that $u$ and $v$ are assumed to be nonconstant, hence $u_M$ and $v_M$ are also nonconstant for all $M$ big enough (in order to use Step 3.2). Setting $\phi_u^M:=\max\{0,u-M\}$, and analogously for $v$, we have
$$
  u=u_M+\phi_u^M,\qquad v=v_M+\phi_v^M.
$$

We will prove below that $E_+[u_M,v_M]=E_+[(u_M)^*,(v_M)^*]$ using the first identity of 
 \eqref{eq:gy:E plus u v is and E minus}. 
Considering $J_-(-t)$ and interchanging the roles of $u$ and $v$, it is easily seen that the same argument will prove the corresponding statement for $J_-$, that is to say, that the second identity of \eqref{eq:gy:E plus u v is and E minus} yields $E_-[u_M,v_M]=E_-[(u_M)^*,(v_M)^*].$
Then, assuming for the moment that the identities $E_\pm[u_M,v_M]=E_\pm[(u_M)^*,(v_M)^*]$ have been proven,
they yield that $E[u_M,v_M]=E[(u_M)^*,(v_M)^*].$ At this point we are in position to apply Step 3.2, which proves that, for every $M>0$, $u_M(x)=(u_M)^{*\eper}(x+z)$ and  $v_M(x)=(v_M)^{*\eper}(x+z)$ for almost every $x\in\R$, for some $z\in\R$ depending a priori on $M$. The fact that $z$ is indeed independent of $M$ follows by an argument as the one below \eqref{eq:uS.uSper.tau}. By letting $M\uparrow+\infty$ and using that $(u^*)_M=(u_M)^{*}$ in $(-\pi,\pi)$ (and analogously for $v$), we finally get that $u(x)=u^{*\eper}(x+z)$ and $v(x)=v^{*\eper}(x+z)$ for some $z\in\R$ and almost every $x\in\R$.

Therefore, it only remains to prove that $E_+[u_M,v_M]=E_+[(u_M)^*,(v_M)^*].$
We first claim that for every $x,y\in(-\pi,\pi)$ the following identity holds:
\begin{equation}
 \label{eq:gy:J ux vy with FM}
 \begin{split}
  J_+(u(x)-v(y))&=J_+(u_M(x)-v_M(y))
  +
  J_+(\phi_u^M(x)-\phi_v^M(y))
   +F_M(\phi_u^M(x),v_M(y)),
 \end{split}
\end{equation}
where
$$
  F_M(\phi,v):=J_+(\phi+M-v)-J_+(\phi)-J_+(M-v).
$$
We will show \eqref{eq:gy:J ux vy with FM} by considering four different cases, depending on where $(x,y)$ lies.
\begin{itemize}
\item \textit{Case} $(x,y)\in \{u>M\}\times \{v>M\}:$ 

In this case 
$$
  u_M(x)=M,\quad v_M(y)=M,\quad u(x)=M+\phi_u^M(x),\quad v(y)=M+\phi_v^M(y),
$$
and 
\begin{equation}
 \begin{split}
  F_M(\phi_u^M(x),v_M(y))=J_+(\phi_u^M(x)+M-M)-J_+(\phi_u^M(x))-J_+(M-M)=0,
\end{split}
\end{equation}
since, as showed in Step 1, we can assume without loss of generality that $J(0)=0$.
From this we get \eqref{eq:gy:J ux vy with FM}.
\smallskip

\item \textit{Case} $(x,y)\in \{u\leq M\}\times\{v>M\}:$ 

In this case $u(x)-v(y)\leq 0 $ and hence $J_+(u(x)-v(y))=0.$ Similarly, as $v_M(y)=M,$ we have $J_+(u_M(x)-v_M(y))=J_+(u(x)-M)=0$, and since $\phi_u^M(x)=0\leq\phi_v^M(y)$ we also get that $J_+(\phi_u^M(x)-\phi_v^M(y))=0.$ Finally, we obtain
$$
  F_M(\phi_u^M(x),v_M(y))=J_+(M-v_M(y))-J_+(0)-J_+(M-v_M(y))=0.
$$

\item \textit{Case} $(x,y)\in \{u>M\}\times\{v\leq M\}:$ 

In this case
$$
  u(x)=M+\phi_u^M(x),\quad v(y)=v_M(y),\quad u_M(x)=M,\quad \phi_v^M(y)=0.
$$
Therefore,
\begin{equation}
\begin{split}
  J_+(u(x)-v(y))&= J_+(M+\phi_u^M(x)-v_M(y)),\\ 
  J_+(u_M(x)-v_M(y))&=J_+(M-v_M(y)),
   \\
  J_+(\phi_u^M(x)-\phi_v^M(y))&=J_+(\phi_u^M(x)),
   \\
  F_M(\phi_u^M(x),v_M(y))&= J_+(\phi_u^M(x)+M-v_M(y))-J_+(\phi_u^M(x))
   -J_+(M-v_M(y)),
  \end{split}
\end{equation}
which proves \eqref{eq:gy:J ux vy with FM} in the present case.
\smallskip

\item \textit{Case} $(x,y)\in \{u\leq M\}\times\{v\leq M\}:$ 

In this case 
$$
  u(x)=u_M(x),\quad v(y)=v_M(y),\quad \phi_u^M(x)=0,\quad \phi_v^M(y)=0,
$$
and
$$
  F_M(\phi_u^M(x),v_M(y))=J_+(M-v_M(y))-J_+(M-v_M(y))=0.
$$
\end{itemize}
In conclusion, \eqref{eq:gy:J ux vy with FM} holds for all $x,y\in(-\pi,\pi)$.

From \eqref{eq:gy:J ux vy with FM} we obtain the following representation of $E_+[u,v]$:
\begin{equation}
 \label{eq:gy:E u v is E uM FM}
  \begin{split}
  E_+[u,v]&=E_+[u_M,v_M]+E_+[\phi_u^M,\phi_v^M] 
  +\int_{-\pi}^{\pi}dx\int_{-\pi}^{\pi}dy\,
  F_M(\phi_u^M(x),v_M(y))g(x-y).
\end{split}
\end{equation}
Note that the third term on the right hand side of \eqref{eq:gy:E u v is E uM FM} involving $F_M$ is nonnegative. This follows from the definition of $F_M$ and the fact that any convex function $f$ with $f(0)=0$ (in our case 
$f:=J_+$) satisfies, for every $a,b\geq 0$ (in our case $a:=\phi_u^M(x)$ and $b:=M-v_M(y)$), that
$$
  f(a+b)-f(a)-f(b)=\int_{a}^{a+b}dt\,f'(t)-f(b)\geq \int_a^{a+b}dt\,f'(t-a)-f(b)
  =0.
$$
Hence all three terms on the right side of \eqref{eq:gy:E u v is E uM FM} are nonnegative and finite.

We will show that the double integral involving $F_M$ on the right hand side of \eqref{eq:gy:E u v is E uM FM} does not increase under the rearrangement of $u$ and $v$. Once we have shown this, we can combine the inequalities $E_+[u_M,v_M]\geq E_+[(u^*)_M,(v^*)_M]$ and  $E_+[\phi_u^M,\phi_v^M]\geq E_+[\phi_{u^*}^M,\phi_{v^*}^M]$ from Step 2.2 (recall that 
$(u_M)^*=(u^*)_M$ and observe that $(\phi_u^M)^*=\phi_{u^*}^M$) with \eqref{eq:gy:E u v is E uM FM} also for $u^*$ and $ v^*$  to obtain, in view of the assumption $E_+[u,v]= E_+[u^*,v^*]$, that all inequalities must be equalities. In particular, we deduce that $E_+[u_M,v_M]= E_+[(u_M)^*,(v_M)^*]$, as desired. 

Thus, it only remains to show that for any pair of functions $\phi\geq 0$ and $0\leq v\leq M$, it holds that
\begin{equation}
 \label{eq:gy:int int FM phi v decr}
  \int_{-\pi}^{\pi}dx\int_{-\pi}^{\pi}dy\,F_M(\phi(x),v(y))g(x-y)
  \geq 
  \int_{-\pi}^{\pi}dx\int_{-\pi}^{\pi}dy\,F_M(\phi^*(x),v^*(y))g(x-y)
\end{equation}
whenever the double integral on the left hand side is finite.
To prove this inequality, let us first observe that the integral of the third term in the definition of $F_M$ is finite, since $J_+$ is nondecreasing in $\R$ and, thus,
\begin{equation}
  \int_{-\pi}^{\pi}dx \int_{-\pi}^{\pi}dy\,J_+(M-v(y))g(x-y)
  \leq 2\pi \|g\|_{L^1(-\pi,\pi)} J(M)<+\infty.
\end{equation}
Therefore, we can  split the left hand side of \eqref{eq:gy:int int FM phi v decr} as 
\begin{equation}
\begin{split}
  \intp dx\intp dy\, F_M(\phi(x)&,v(y))g(x-y)
  =  \Sigma_M[\phi,v]-\int_{-\pi}^{\pi}dx\int_{-\pi}^{\pi}dy\,J_+(M-v(y))g(x-y),
\end{split}
\end{equation}
where
$$
  \Sigma_M[\phi,v]:=\int_{-\pi}^{\pi}dx\int_{-\pi}^{\pi}dy\,
  \Big(J_+(\phi(x)+M-v(y))-J_+(\phi(x))\Big)g(x-y).
$$
To prove \eqref{eq:gy:int int FM phi v decr}, it is sufficient to show that $\Sigma_M$ does not increase under rearrangement because, using \eqref{eq:gy:F abs contin ustar} with $F(t):=J_+(M-t)$, we have
\begin{equation}
  \int_{-\pi}^{\pi}dy\,J_+(M-v(y))
  =\int_{-\pi}^{\pi}dy\,J_+(M-|v(y)|)=\int_{-\pi}^{\pi}dy\,J_+(M-v^*(y)).
  \end{equation}
At this point, we now proceed in a similar way as how we obtained \eqref{eq:gy:E u v finite and Fubini}. More precisely, we use the identity
\begin{equation}
\begin{split}
  J_+(\phi(x)+M-v(y))-J_+(\phi(x))&=
  \int_{\phi(x)}^{\phi(x)+M-v(y)}d\tau\,J_+'(\tau)
  =\int_0^{M-v(y)}d\tau\,J_+'(\phi(x)+\tau)\\
  &=\int_0^Md\tau\,J_+'(\phi(x)+\tau)(1-\chi_{\{v>M-\tau\}}(y))
\end{split}
\end{equation}
to write
$$
  \Sigma_M[\phi,v]=\int_0^M d\tau\int_{-\pi}^{\pi}dx\int_{-\pi}^{\pi}dy\,J_+'(\phi(x)+\tau)g(x-y)(1-\chi_{\{v>M-\tau\}}(y)).
$$
Finally, that $\Sigma_M$ does not increase under rearrangement follows from \Cref{thm:gy:sharp Friedberg Luttinger} and the monotone convergence theorem---exactly as in Step 2, where we have shown the same result for $E_+$ instead of $\Sigma_M$. In conclusion, \eqref{eq:gy:int int FM phi v decr} holds and Step 3.3 is now complete.

\medskip
\noindent\underline{{\em Step 4} (case of equality for $J(t)=|t|$):} 
\smallskip

Since $J_+'(u(x)-\tau)=1$ if $u(x)> \tau$ and zero elsewhere, using \eqref{eq:gy:A and B} and \eqref{eq:gy:E u v split in A B and ineq} we see that $E[u,v]=E[u^*,v^*]$ if and only if,  for almost every $\tau>0$,
$$
  \intp dx\intp dy\, \chi_{\{u> \tau\}}(x)g(x-y)\chi_{\{v>\tau\}}(y)
  =
  \intp dx\intp dy\, \chi_{\{u^*> \tau\}}(x) g(x-y)\chi_{\{v^*>\tau\}}(y).
$$
From this it follows immediately that \eqref{rear:|t|.sets} yields $E[u,v]=E[u^*,v^*]$. So it remains to show the converse.

Note that neither of the two functions $\chi_{\{u> \tau\}}$ and $\chi_{\{v>\tau\}}$ is constant almost everywhere if and only if
\begin{equation}\label{eq:gy:tau in ess inf ess sup}
  \tau\in (\operatorname{ess}\inf u,\operatorname{ess}\sup u)\cap (\operatorname{ess}\inf v,\operatorname{ess}\sup v).
\end{equation}
Thus, from the case of equality in  \Cref{thm:gy:sharp Friedberg Luttinger} we obtain that, 
for almost every $\tau$ in the range given by \eqref{eq:gy:tau in ess inf ess sup} there exists $z_{\tau}\in\R$ such that 
\begin{equation}\label{eq:case_eq_|t|_rangetau}
  \chi_{\{u> \tau\}}(x)=\chi_{\{u^{*\eper}> \tau\}}(x+z_{\tau}),\qquad \chi_{\{v>\tau\}}(x)=\chi_{\{v^{*\eper}>\tau\}}(x+z_{\tau})
  \quad\text{ for a.e.\! }x\in\R.
\end{equation}
In particular, for almost every $\tau$ as in \eqref{eq:gy:tau in ess inf ess sup}, $\{u> \tau\}$ and $\{v> \tau\}$ are ``$2\pi$-periodic unions of intervals''. Then, thanks to  \Cref{lemma:gy:slicing of layer cake}~$(ii)$ applied to each of these intervals, \eqref{eq:case_eq_|t|_rangetau} actually holds for every $\tau$ as in \eqref{eq:gy:tau in ess inf ess sup}.

\medskip

\noindent\underline{{\em Step 5} (case that 
$\inf J$ is not attained):} 
\label{Step inf J not att}
\smallskip

Throughout Step 1 to 4 we were assuming that $\inf J$ was attained, that is, that there exists $t_0\in\R$ such that 
$J(t_0)\leq J(t)$ for all $t\in\R$, and we proved \Cref{theorem:gy:non expansivity under J} in this particular case. It remains to prove the theorem in the case when $\inf J$ is not attained. Since $J$ is nonnegative by assumption,  we have that $0\leq \inf J<+\infty$ and thus, without loss of generality, we can assume that $\inf J=0,$ arguing as in \eqref{eq:minJ_unnecessary} in Step 1. Note also that a convex function in the real line which does not attain its infimum must be either increasing or decreasing (this can be seen arguing by contradiction). Therefore, we can also assume without loss of generality that $J$ is increasing and, based on the previous observation, that $J(-\infty):=\lim_{t\downarrow-\infty}J(t)=0$; in the case that $J$ is decreasing simply consider $J(-t)$. Now the splitting $J=J_+ +J_-$ is not necessary (and not possible), but since $J(-\infty)=0$ we have as in Step 2 that
\begin{equation}
\begin{split}
 J(u(x)-v(y))&= \int_{-\infty}^{u(x)-v(y)}d\tau\, J'(\tau)=
  \int_{v(y)}^{+\infty}d\tau\, J'(u(x)-\tau)\\ 
  &=
  \int_{0}^{+\infty}d\tau\, J'(u(x)-\tau)\left(1-\chi_{\{v>\tau\}}(y)\right).
\end{split}
\end{equation}
Thus, we obtain the same expression for $E[u,v,g]$ as in \eqref{eq:gy:E u v finite and Fubini}, with the only difference that $J_+$ is replaced by $J$, but recall that now $J'$ is nondecreasing (by the convexity of $J$) and nonnegative (since $J$ is increasing) in $\R$, as it was $J'_+$ in Step 2. With this in hand, the inequality \eqref{rear:Euv} in \Cref{theorem:gy:non expansivity under J} can be proven exactly as in Step 2. 

Next, we study the case of equality. First observe that in Step 3.2 we did not use the assumption that $\inf J$ is attained, nor the splitting  $J=J_++J_-$. Hence we only have to deal with the case of equality in the case that $u$ or $v$ (or both) are unbounded. Let us first assume that only $v$ might be unbounded, but that $u$ is bounded. We will use the same abbreviations and argument as in Step 3.3. Let us write
$$
  J(u(x)-v(y))=J(u(x)-v_M(y))+J(u(x)-M-\phi_v^M(y))-J(u(x)-M).
$$ 
It is immediate to check this identity by distinguishing the two cases $v(y)\leq M$ and $v(y)>M.$
If we define this time
$$
  \Sigma_M[u,\phi]:=\int_{-\pi}^{\pi}dx\int_{-\pi}^{\pi}dy\,
  J(u(x)-M-\phi(y))g(x-y),
$$
then we obtain
$$
  E[u,v]=E[u,v_M]+\Sigma_M[u,\phi_v^M]-\int_{-\pi}^{\pi}dx\int_{-\pi}^{\pi}dy\, J(u(x)-M)g(x-y).
$$
The third integral is finite, because $J$ is increasing and therefore $J(u(x)-M)\leq J(\esssup u-M)$ for almost every $x$. By assumption $E[u,v]$ is finite, and thus all three terms on the right hand side of the previous expression are finite too (because $E$ and $\Sigma_M$ are nonnegative).
As in Step 3.3, using \eqref{eq:gy:F abs contin ustar} we see that the third term involving $J(u(x)-M)$ on the right hand side is unaltered under the rearrangement of $u.$ Moreover, 
$E[u,v_M]$ and $\Sigma_M[u,\phi_v^M]$ do not increase under the rearrangements of $u$, $v_M$, and $\phi_v^M$, in view of the inequality of the theorem. Hence the equality $E[u,v]=E[u^*,v^*]$ forces the equality $E[u,v_M]=E[u^*,(v_M)^*]$ to hold. Now, $u$ and $v_M$ are both bounded and we can conclude as in Step 3.2, which, we recall, does not use the assumption that $\inf J=\min J$.

Let us finally assume that both $u$ and $v$ are unbounded, which is the last remaining case. We will repeat the previous argument, using that \Cref{theorem:gy:non expansivity under J} holds true for general $v,$ but $u$ bounded. This time use  the identity
\begin{align}
  J(u(x)-v(y))=J(u_M(x)-v(y))+J(\phi_u^M(x)+M-v(y))-J(M-v(y)),
\end{align}
and set
$$
  \Sigma_M[\phi,v]:=\int_{-\pi}^{\pi}dx\int_{-\pi}^{\pi}dy\,
  J(\phi(x)+M-v(y))g(x-y).
$$
Note that  $J(M-v(y))\leq J(M)<+\infty$ because $J$ is increasing. Thus, we can proceed exactly as in the previous case and show that the equality $E[u,v]=E[u^*,v^*]$ forces the equality $E[u_M,v]=E[(u_M)^*,v^*]$ to hold. Since we have already dealt with the case that one of the two functions is bounded, $u_M$ in the present case, we are done letting $M\uparrow+\infty$ as in the beginig of Step 3.3.  
\medskip

\end{proof}
\smallskip

\section{Nonexpansivity under Schwarz rearrangement}
\label{section nonexpansive cyl}

We will now prove the nonexpansivity result under standard Schwarz symmetrization, namely \Cref{theorem:gy:J in R n-1}. As already mentioned in the introduction, we will only detail those steps which differ from the proof of the periodic case given in  \Cref{section nonexpansive periodic}. Let us use the abbreviations
\begin{equation}
\begin{split}
  \overline{E}[u,v,g]&:=\int_{\R^{n}}dx\int_{\R^{n}}dy\,J(u(x)-v(y))g(x-y),
  \\
  \overline{E}[u,v]&:=\overline{E}[u,v,g]\quad\text{ if $g=g^{*\,n}$},
 \end{split}
\end{equation}
for two nonnegative measurable functions $u,v:\R^{n}\to\R.$ Before giving the proof, we shall mention some elementary observations on the hypothesis of the theorem for the sake of completeness.

Note that the assumption $|\{u>\tau\}|<+\infty$ is necessary to give sense to the theorem (in contrary to the periodic case), otherwise the rearrangement is in general not defined.

Concerning the hypothesis of nonnegativeness of $u$ and $v$, a similar comment as \Cref{remark:negative functions} for the periodic version of the theorem is valid here. For instance if $J(t)=|t|^p,$ $p\geq 1,$ then the inequality remains true even for sign-changing functions $u$ and $v.$ But the hypothesis of nonnegativeness cannot be relaxed for the statement $(i)$ of the theorem to hold true. To see this we give a different example than the one of \Cref{remark:negative functions}, which does not work in the present setting, since infinite integrals would appear if we take $u$ or $v$ to be a nonzero constant function. Instead, consider the strictly convex function $J$ given by
$$
  J(t)=\left\{\begin{array}{rl}
               t^2 & \text{ if }t\leq 0,
               \\
               \frac{t^2}{2} &\text{ if }t> 0.
              \end{array}\right.
$$
Let $v\equiv 0$ and $u$ be some nonpositive function with compact support, not identical to zero. Then,
$$
  \int_{\R^n}J(u^{*\,n})=\int_{\R^n} \frac{(u^{*\,n})^2}{2}=\int_{\R^n}\frac{|u|^2}{2}<\int_{\R^n} |u|^2=\int_{\R^n}J(u),  
$$
which yields $\overline{E}[u^{*\,n},v^{*\,n}]<\overline{E}[u,v].$ Thus, $\overline{E}$ decreases under the rearrangement  although $v$ is identically zero.

The proof of \Cref{theorem:gy:J in R n-1} will actually show that the conclusion of $(ii)$ holds true if one only requires that $J_+:=\chi_{[0,+\infty)}J$ is strictly convex and the supremum of $u$ is greater than or equal to the supremum of $v$ (or, similarly, that $J_-:=J-J_+$ is strictly convex and $\esssup v\geq \esssup u$).

\begin{proof}[Proof of \Cref{theorem:gy:J in R n-1}]
The proof follows the same outline as that of  \Cref{theorem:gy:non expansivity under J}, hence we will only carry out the steps which are different. The decomposition  $J=J_++J_-$ is identical as in the proof of \Cref{theorem:gy:non expansivity under J}, and also expressing $\overline{E}_+$ in terms of the derivative of $J_+$.
Since $\overline{E}_+[u,v,g]$ is finite, Fubini's theorem gives
$$
  \overline{E}_+[u,v,g]=\int_0^{+\infty}d\tau\,(A(u,g,\tau)-B(u,v,g,\tau)),
$$
where $A$ and $B$ are defined as in \eqref{eq:gy:A and B} with the only difference that the integrals are over $\R^n\times\R^n$ instead of $(-\pi,\pi)\times(-\pi,\pi)$, that is,
\begin{equation}
\begin{split}
  A(u,g,\tau)&:=\int_{\R^n}dx\int_{\R^n}dy\,
  J'_+(u(x)-\tau)g(x-y),
  \smallskip \\
  B(u,v,g,\tau)&:=\int_{\R^n}dx\int_{\R^n}dy\,
  J'_+(u(x)-\tau)g(x-y)\chi_{\{v>\tau\}}(y).
  \end{split}
\end{equation}

\medskip
\noindent\underline{{\em Step 1} (proof of the inequality):} 
\smallskip

Let us first assume that either $u$ is bounded or that $|J'|$ is bounded ---the latter occurring for instance if $J(t)=|t|.$ 
Note that if $\esssup u=M<+\infty$
then $A(u,g,\tau)$ is finite because $g\in L^1(\R^n)$ and
$$
  \int_{\R^n}dx\,J'_+(u(x)-\tau)=\int_{\{u>\tau\}}dx\,J'_+(u(x)-\tau)
  \leq |\{u>\tau\}| J_+'(M-\tau)<+\infty.
$$
Similarly, if $|J'|$ is bounded one obtains as well that $A(u,g,\tau)$ is finite. As in the proof of  \Cref{theorem:gy:non expansivity under J}, using  \Cref{lemma:gy:increasing and low semicont}, we get $A(u,g,\tau)=A(u^{*\,n},g^{*\,n},\tau).$ Since $A(u,g,\tau)$ is finite we get that $B(u,v,g,\tau)$ is finite too for almost every $\tau$, and therefore  
\begin{equation}
 \label{eq:gy:overline E and E ast difference}
  \overline{E}_+[u,v,g]=\overline{E}_+[u^{*\,n},v^{*\,n},g^{*\,n}]
  +\int_0^{+\infty}d\tau\,(B(u^{*\,n},v^{*\,n},g^{*\,n},\tau)-B(u,v,g,\tau)),
\end{equation}
where $B(u^{*\,n},v^{*\,n},g^{*\,n},\tau)$ and $B(u,v,g,\tau)$ are finite for almost every $\tau>0.$
Hence it is sufficient to show that $B(u,v,g,\tau)\leq B(u^{*\,n},v^{*\,n},g^{*\,n},\tau).$ This last inequality follows from the classical Riesz rearrangement inequality \eqref{eq:gy:intro:classical Riesz rearrangement}; see for instance \cite[Theorem 3.7]{Lieb Loss}.  This proves \eqref{cyl.rear.ineq}  in the case that $u$ or $|J'|$ is bounded. 

If neither $u$ nor $|J'|$ is bounded, then one can proceed exactly as in Step 2.2 of the proof of  \Cref{theorem:gy:non expansivity under J} by defining $u_M$ as in \eqref{eq:gy:u M Def} and using the monotone convergence theorem.

\medskip
\noindent\underline{{\em Step 2} (case of equality for strictly convex $J$ under a boundedness assumption):} 
\smallskip

From now on, we will omit in all the functionals $\overline{E},A,B,\ldots$ the argument $g=g^{*\,n},$ as it is unaltered under rearrangement.

If $u$ is constant, then by the hypothesis $|\{u>\tau\}|<+\infty$ for every $\tau>0$ we must have that $u\equiv 0$. In this case there is equality in \eqref{cyl.rear.ineq} for any $v$ (whether $J$ is strictly convex or not). Since the same argument applies if $v$ is constant, from now on we can assume that neither $u$ nor $v$ is constant. 

In this step we will study the case of equality under the assumption that either $u$ and $v$ are both bounded or that $|J'|$ is bounded. By Step 1 it holds that 
\begin{equation}
 \label{eq:gy:E u v split in A B and ineq Rm}
  \overline{E}_+[u,v]=\overline{E}_+[u^{*\,n},v^{*\,n}]+\int_0^{+\infty}d\tau\,
  (B(u^{*\,n},v^{*\,n},\tau)-B(u,v,\tau))
\end{equation}
and, therefore, that $B(u,v,\tau)=B(u^{*\,n},v^{*\,n},\tau)$ for almost every $\tau>0.$ Note that, from the assumption that $|\{u>\tau\}|$ and $|\{v>\tau\}|$ are finite for all $\tau>0,$ we see that $\essinf u=\essinf v=0.$ Let us now assume $\esssup u\geq\esssup v $ (for the other case, in all of the following argument interchange the roles of $J_+(t)$ and $J_-(-t)$). Hence we have
$$
  \mathcal{R}(v):=(\essinf v,\esssup v)\subset(\essinf u,\esssup u)
  =:\mathcal{R}(u).
$$

We now apply the strict Riesz rearrangement inequality to the second term on the right hand side of \eqref{eq:gy:E u v split in A B and ineq Rm}; see Baernstein \cite[Theorem 2.15 $(b)$]{Baernstein book} or Lieb and Loss \cite[Theorem 3.9]{Lieb Loss}.\footnote{\cite[Theorem 3.9]{Lieb Loss} should include two hypotheses which are needed for its validity —and which are explicitly
stated in most of the results of \cite[Chapter 3]{Lieb Loss}. Namely, that the functions $f$ and $h$ of its statement must vanish
at infinity and that neither $f$ nor $h$ is identically zero. Indeed, there is always equality in the Riesz rearrangement
inequality if one of these two functions is constant.} 
Observe first that if $\tau\in \mathcal{R}(v)$ then neither $J'_+(u-\tau)$ nor $\chi_{\{v>\tau\}}$ is constant (in particular, they are not identical to zero), and we therefore obtain from $B(u,v,\tau)=B(u^{*\,n},v^{*\,n},\tau)$ that there exists $a_{\tau}\in\R^{n}$ depending a priori on $\tau$ such that
\begin{equation}
 \label{eq:gy:J prime plus and chi v}
  J_+'(u(x)-\tau)=J_+'(u^{*\,n}(x-a_{\tau})-\tau)
  \,\text{ and }\,
  \chi_{\{v>\tau\}}(y)=\chi_{\{v^{*\,n}>\tau\}}(y-a_{\tau})
\end{equation}
for a.e.\! $x\in\R^n$ and for a.e.\! $y\in\R^n$.
 Since $J_+$ is strictly convex in $[0,+\infty)$, we have that $J_+'$ is increasing and, therefore, one-to-one in $[0,+\infty)$. Using this together with the fact that $J'_+$ is identical to zero in $(-\infty,0)$, we deduce that, for almost every $\tau\in\mathcal{R}(v)$,
\begin{equation}
 \label{eq:gy:u geq tau is u ast}
  \{u>\tau\}=\{u^{*\,n}>\tau\}+a_{\tau} \quad\text{in measure}
\end{equation}
and
\begin{equation}
 \label{eq:gy:u ast m minut tau}
  u(x)=u^{*\,n}(x-a_{\tau})\quad \text{for a.e.\! }x\in \{u^{*\,n}>\tau\}+a_{\tau}.
\end{equation}
We will now show that $a_{\tau}$ is indeed independent of $\tau$ and thus, taking $\tau\downarrow0$ in \eqref{eq:gy:u ast m minut tau}, that $u$ is a translate of a symmetric decreasing function.\footnote{Let us point out the difference with Step 3.2 of the proof of the corresponding theorem in the periodic case. There we had in \eqref{eq:uS minus tau} that $u_S(s)-\tau>0$ for all $\tau$ in the given range  and \eqref{eq:gy:u ast m minut tau} followed for all $x$. This trick of ``shifting'' $u$ by defining $I_S$ does not work in the case where $(-\pi,\pi)$ is replaced by $\R^n,$ because some integrals might not be finite now. However, this was not the only reason why we introduced $I_S$: it was also due to the fact that in the periodic setting we may not have $\essinf u=\essinf v=0$, which causes some technical difficulties.}

Suppose that $s,\tau\in \mathcal{R}(v)$ are such that $s<\tau$. We claim that \eqref{eq:gy:u geq tau is u ast} and \eqref{eq:gy:u ast m minut tau} applied to $s$ and $\tau$ gives
\begin{equation}
 \label{eq:gy:u geq tau shift by as}
  \{u>\tau\}=\{u^{*\,n}>\tau\}+a_s\, \quad\text{in measure}.
\end{equation}
Let us first show that, up to a set of measure zero, the set on the left hand side is included in the set on the right hand side. If $x\in\{u>\tau\}\subset\{u>s\}$ then by \eqref{eq:gy:u geq tau is u ast} $x\in\{u^{*\,n}>s\}+a_s$. We now use \eqref{eq:gy:u ast m minut tau} with $s$ and obtain that $\tau<u(x)=u^{*\,n}(x-a_s),$ and thus $x\in \{u^{*\,n}>\tau\}+a_s$. This proves that
$$
  \{u>\tau\}\subset\{u^{*\,n}>\tau\}+a_s\,.
$$
Now the two sets have the same measure, and thus they have to be equal in measure, which proves \eqref{eq:gy:u geq tau shift by as}. Since the sets $\{u^{*\,n}>\tau\}$ are nonempty bounded balls centered at the origin for every $\tau\in \mathcal{R}(v)\subset\mathcal{R}(u)$, the identities \eqref{eq:gy:u geq tau is u ast} and \eqref{eq:gy:u geq tau shift by as} force $a_{\tau}=a_s$. This shows that $a_{\tau}\equiv a$ is independent of $\tau,$ which gives, by letting $\tau\downarrow 0=\essinf u$ and using \Cref{lemma:gy:slicing of layer cake}~$(i)$, that
$$
  u(x)=u^{*\,n}(x-a)\quad\text{for a.e.\! }x\in \R^n.
$$

Finally, we now use the second identity in \eqref{eq:gy:J prime plus and chi v} to obtain that 
$$
  \chi_{\{v>\tau\}}(y)=\chi_{\{v^{*\,n}>\tau\}}(y-a)\quad\text{for a.e.\! }y\in \R^n.
$$
Since this holds for almost every $\tau\in \mathcal{R}(v)$, we obtain from  \Cref{lemma:gy:slicing of layer cake} $(i)$  that $v(y)=v^{*\,n}(y-a)$ for a.e.\! $y\in \R^n$.

\medskip
\noindent\underline{{\em Step 3} (case of equality for $J(t)=|t|$):} 
\smallskip

The proof of this case is exactly the same as that of the corresponding result in  \Cref{theorem:gy:non expansivity under J}. Observe that since $|\{u>\tau\}|$ and $|\{v>\tau\}|$ are finite for all $\tau>0$ by assumption, we have 
$\operatorname{ess}\inf u=\operatorname{ess}\inf v=0.$

\medskip
\noindent\underline{{\em Step 4} (case of equality for unbounded $u$ or $v$, and unbounded $|J'|$)}\footnote{We believe this proof is a simplification of the proof given in \cite[Lemma A.2]{Frank Seiringer} in the special case $u=v$, where the second derivative of $J$ was used to carry out this step.} 
\smallskip 

In Step 2 we studied the case of equality under the assumption that either $u$ and $v$ are both bounded or that $|J'|$ is bounded. The purpose of this step is study the case of equality without this assumption.

  Let $t_i$ be a sequence of positive real numbers such that $\lim_{i\uparrow+\infty}t_i=+\infty$ and $J_+$ is differentiable at $t_i$ (the use of points of differentiability is not really necessary to carry out this step but it makes the proof less heavy). Define $J_{+,i}$ as the convex function in $\R$ which coincides with $J_+$ in $(-\infty,t_i]$ and is a linear function with slope $J'_+(t_i)$  in $(t_i,+\infty)$, namely
$$
  J_{+,i}(t):=J_+(t)\quad\text{if }t\leq t_i\,,\quad J_{+,i}(t):=J_+(t_i)+(t-t_i)J'_+(t_i)\quad\text{if }t> t_i\,.
$$
This gives an increasing sequence $J_{+,i}(t)\leq J_{+,i+1}(t)\leq J_+(t)$ wich converges to $J_+(t)$ for every $t\in\R$ as $i\uparrow+\infty$. A analogous cutt-off $J_{-,i}$ can be done for $J_-,$ choosing some sequence of negative real numbers $s_i\downarrow-\infty$ and replacing $J_-$ in $(-\infty,s_i)$ by a suitable linear function.
Let us define
$$
  \overline{E}^i_{\pm}[u,v]:=\int_{\R^n}dx\int_{\R^n}dy\,J_{\pm,i}(u(x)-v(y))g(x-y).
$$
We claim that if $\overline{E}$ is unaltered under the rearrangement of $u$ and $v$, then for every $i$ we have
\begin{equation}
 \label{eq:gy:overline Ei pm equalities}
  \overline{E}_+^i[u,v]=\overline{E}_+^i[u^{*\,n},v^{*\,n}]
  \quad\text{ and }\quad
  \overline{E}_-^i[u,v]=\overline{E}_-^i[u^{*\,n},v^{*\,n}].
\end{equation}
First observe that the equality $\overline{E}[u,v]=\overline{E}[u^{*\,n},v^{*\,n}]$ yields the validity if $\eqref{eq:gy:overline Ei pm equalities}$ without the superscript $i$ because $\overline{E}_\pm$ do not increase under Schwarz rearrangement. Hence, we only need to show that from the equality $\overline{E}_+[u,v]=\overline{E}_+[u^{*\,n},v^{*\,n}]$ it follows the first equality of \eqref{eq:gy:overline Ei pm equalities},  since the corresponding second identity of \eqref{eq:gy:overline Ei pm equalities} would then follow by considering $J(-t).$ 

Let us decompose $\overline{E}_+=\overline{E}_+^i+\overline{\Sigma}_i$, where 
$$
  J_+=J_{+,i}+\Phi_i \quad\text{and}\quad
  \overline{\Sigma}_i[u,v]:=\int_{\R^n}dx\int_{\R^n}dy\,\Phi_i(u(x)-v(y))g(x-y).
$$
Observe that $\Phi$ is a convex function, because it is identically zero in $(-\infty,t_i]$ and in $[t_i,+\infty)$ it is the (nonnegative) difference of a convex function and a linear function. Hence, by Step 1, $\overline\Sigma_i$ does not increase under rearrangement.
Now, by assumption we have that 
$$
  \overline{E}_+^i[u,v]+\overline{\Sigma}_i[u,v]
  =\overline{E}_+^i[u^{*\,n},v^{*\,n}]  +\overline{\Sigma}_i[u^{*\,n},v^{*\,n}],
$$
where all terms are nonnegative. Since both $\overline{E}_+^i$ and $\overline{\Sigma}_i$ do not increase under rearrangement, we obtain \eqref{eq:gy:overline Ei pm equalities} for $\overline{E}_+^i$, as desired.

Since $|J_{\pm,i}'|$ are bounded we can now apply the argument of Step 2 to the identities of \eqref{eq:gy:overline Ei pm equalities}. The first identity of \eqref{eq:gy:overline Ei pm equalities} gives for instance that,  for a.e.\! $\tau\geq 0$,
$$
  \int_{\R^n}dx\int_{\R^n}dy\,
  J'_{+,i}(u(x)-\tau)g(x-y)\chi_{\{v>\tau\}}(y)
  =
  \int_{\R^n}dx\int_{\R^n}dy\,
  J'_{+,i}(u^{*\,n}(x)-\tau)g(x-y)\chi_{\{v^{*\,n}>\tau\}}(y).
$$
Let us assume as in Step 2 that the supremum of $u$ is greater than or equal to the supremum of $v$, thus $\operatorname{ess}\sup u=+\infty$ because $u$ or $v$ is unbounded by assumption.  Observe that  $J_{+,i}$ is strictly convex only in $(0,t_i).$ Therefore the arguments of Step 2 have to be actually slightly modified (this is the only difference to Step 2). If we define
$
  \mathcal{R}(v):=(0,\operatorname{ess}\sup v),
$
then the analogue of \eqref{eq:gy:J prime plus and chi v} holds for almost every $\tau\in \mathcal{R}(v)$ ---since 
$\tau<\esssup v$ we have that $\chi_{\{v>\tau\}}$ is not indentical to zero, and since 
$\essinf u=0$, $\esssup u=+\infty$, and $J'_{+,i}(t)=J'_{+}(t_i)>0$ for all $t>t_i$ we have that $J_{+,i}'(u-\tau)$ is not indentical to zero--- and rewrites as
\begin{equation}
  J_{+,i}'(u(x)-\tau)=J_{+,i}'(u^{*\,n}(x-a_{\tau,i})-\tau)
  \,\text{ and }\,
  \chi_{\{v>\tau\}}(y)=\chi_{\{v^{*\,n}>\tau\}}(y-a_{\tau,i}).
\end{equation}
From the first identity one deduces (since $J_{+,i}'$ is one-to-one in $(0,t_i)$ and if $u<t_i$ then $u-\tau<t_i$) that, as in \eqref{eq:gy:u geq tau is u ast},
$$
  \{\tau<u<t_i\}=\{\tau<u^{*\,n}<t_i\}+a_{\tau,i}
  \quad\text{ in measure,}
$$
and
\begin{equation}\label{lim.tau.u.t.i}
  u(x)=u^{*\,n}(x-a_{\tau,i})\quad\text{for a.e.\! }x\in 
  \{\tau<u^{*\,n}<t_i\}+a_{\tau,i},
\end{equation}
for a.e.\! $\tau\in \mathcal{R}(v)$.
Arguing as in \eqref{eq:gy:u geq tau shift by as} one  obtains that if $0<s\leq \tau$ then
$$
  \{\tau<u^{*\,n}<t_i\}+a_{\tau,i}=\{\tau<u^{*\,n}<t_i\}+a_{s,i}.
$$
Since $\{\tau<u^{*\,n}<t_i\}$ is either a bounded annulus, or the empty set (if $t_i\leq\tau$ we have an  empty set, and if $t_i>\tau$ we have an annulus since $\esssup u=+\infty$ and $\{u>\tau\}$ has finite measure for all $\tau>0$) we obtain that $a_{\tau,i}\equiv a_i$ is independent of $\tau.$ Letting $\tau\downarrow0$ in \eqref{lim.tau.u.t.i} gives
$$
  u(x)=u^{*\,n}(x-a_i)\quad\text{for a.e.\! $x\in\{u<t_i\}$.}
$$
One easily sees that $a_i\equiv a$ cannot depend on $i$, and by letting $t_i\uparrow+\infty$ we obtain that $u=u^{*\,n}(x-a).$

Exactly as in Step 2, we also obtain that $v(y)=v^{*\,n}(y-a).$
\end{proof}

\section{Rearrangements and periodic Gagliardo seminorms}
\label{subsection Gagliardo seminorms}

Let us define, for $1\leq p<+\infty$ and a function $u:\R^n\to\R$ which is $2 \pi$-periodic in the variable $x_1$, the following periodic Gagliardo seminorm
\begin{equation}
  [u]_{W^{s,p}}^{\eper}:=\left(\int_{\{x\in\R^n:\,-\pi<x_1<\pi\}}dx\int_{\R^n}dy\,\frac{|u(x)-u(y)|^p}{|x-y|^{n+s p}} \right)^{1/p}.
\end{equation}
Indeed, if a set $E\subset\R^n$ is $2 \pi$-periodic in the variable $x_1$, then taking $u$ equal to the characteristic function $\chi_E$ gives that $[\chi_E]_{W^{s,1}}^{\eper}=2\PP_s(E)$,
where  
$$
   \PP_s[E] := \int_{E\cap\{-\pi<x_1<\pi\}}\!dx
   \int_{\R^n\setminus E}dy\, \frac{1}{|x-y|^{n+s}}
$$
is the periodic fractional perimeter functional introduced in \cite{Davila}; 
here we have also used that $\PP_s(\R^n\setminus E)=\PP_s(E)$ by the periodicity of $E$. 
More generally, one can show that for every function $u$ which is $2\pi$-periodic in the variable $x_1$ it holds that
\begin{equation}
  \label{eq:gy:W s1 period and PPs}
  [u]_{W^{s,1}}^{\eper}=2\int_{-\infty}^{+\infty}dt\,\PP_s(\{u>t\}).
\end{equation}
We call this identity the  periodic fractional coarea formula, which follows from the identity
$$
  |u(x)-u(y)|=\int_{-\infty}^{+\infty}dt\left(\chi_{\{u>t\}}(x)\chi_{\{u\leq t\}}(y) +\chi_{\{u\leq t\}}(x)\chi_{\{u>t\}}(y)\right)
$$
and Fubini's theorem.

Our first result of this section deals with the behavior of the  periodic Gagliardo seminorm under periodic rearrangement. Recall that $u^{*\eper}(x_1,x')=(u(\cdot,x'))^{*\eper}(x_1)$ where, for a periodic function $f$ in $\R$, $f^{*\eper}$ is simply the standard Schwarz (or Steiner) symmetrization of $f$ resticted to one period centered at $0$ and then extended in a periodic way.

\begin{theorem}
\label{thm:gy:periodic nonlocal Polya}
Let $n\geq 1$ be an integer, $0<s<1,$ $1\leq p<+\infty$,  and $u:\R^n\to\R$ be a measurable function that is $2\pi$-periodic in the variable $x_1$ and with finite seminorm $[u]_{W^{s,p}}^{\eper}.$ Then, 
\begin{equation}
 \label{eq:periodic nonlocal Polya ineq}
  \left[u^{*\eper}\right]_{W^{s,p}}^{\eper}
  \leq
   \left[u\right]_{W^{s,p}}^{\eper}.
\end{equation}

Moreover, the following holds:

\begin{itemize}
\item[$(i)$] If $p>1,$ or $p=1$ and $u$ is the characteristic function of some measurable set $E,$ then there is equality in  \eqref{eq:periodic nonlocal Polya ineq} if and only if $u=\pm|u|$ and, modulo translations, $|u|$ is nonincreasing in $(0,\pi)$ in the variable $x_1$ and even with respect to $\{x_1=0\}.$ More precisely, there is equality in  \eqref{eq:periodic nonlocal Polya ineq} if and only if
$$
  u(x_1,x')=\pm u^{*\eper}(x_1+a,x')\quad\text{for a.e.\! }(x_1,x')\in\R\times\R^{n-1},
$$
for some $a\in\R.$

\item[$(ii)$] If $p=1$ and $u$ is continuous then there is equality in \eqref{eq:periodic nonlocal Polya ineq} if and only if $u=\pm|u|$ and for almost every $\tau \in\R$ there exists $a_{\tau}\in\R$ such that the superlevel sets satisfy $$\{|u|>\tau\}=\{|u|>\tau\}^{*\eper}+a_{\tau}e_1,$$
up to redefining $u$ on a set of vanishing $n$-dimensional Lebesgue measure. Here $e_1:=(1,0,\ldots,0)$ denotes the first coordinate direction of $\R^n$.
\end{itemize}
\end{theorem}

One of the statements of $(ii)$, namely the 
``if'', 
does not require that $u$ is continuous. Any function which is of the form described in $(ii)$, whether continuous or not, preserves the periodic Gagliardo seminorm 
$[\,\cdot\,]_{W^{s,1}}^{\eper}$ under periodic rearrangement. This follows from $(i)$ and the periodic fractional coarea formula \eqref{eq:gy:W s1 period and PPs}.
Note that if $u$ is as in $(ii)$, this does not force $u$ to be (modulo translations) symmetric and nonincreasing in the variable $x_1$ in $(-\pi,\pi)\times\R^{n-1}$, even if all its superlevel sets are. Consider for example the function $u:(-2,2)\to \R,$ given by $u(x)=\chi_{(-1,2)}(x)+\chi_{(-1,0)}(x)$, and extended in a $4$-periodic way to $\R$. Then $u^{*\eper}$ is, up to a translation,  the $4$-periodic extension of $\chi_{(-1,2)}+\chi_{(0,1)}.$ It is easy to check that for $p=1$ the periodic Gagliardo seminorms of $u$ and $u^{*\eper}$ are equal.
On the other hand, we were not able to establish whether the continuity hypothesis on $u$ is necessary or not in the ``only if'' statement of $(ii)$.

Our second result of this section concerns the behavior of the periodic Gagliardo seminorm under a cylindrical rearrangement, more precisely under the spherical rearrangement in $\{x_1\}\times\R^{n-1}$ for every frozen $x_1\in\R$.  We recall here its definition (given in \Cref{section preliminaries}), which is 
$$ 
  u^{*\,n,1}(x_1,x'):=u(x_1,\cdot)^{*\,n-1}(x'),
$$
where $v^{*\,n-1}$ denotes the usual Schwarz symmetrization (spherical decreasing rearrangement) of $v:\R^{n-1}\to\R.$
In other words, $u^{*\,n,1}$ is cylindrically symmetric around the $x_1$-axis and nonincreasing in $|x'|.$ As stated in the following theorem, this rearrangement also decreases the periodic Gagliardo seminorm, although this has nothing to do with periodicity.

\begin{theorem}
\label{thm:gy:Polya for n minus 1 Schwarz}
Let $n\geq 2$ be an integer, $0<s<1,$ $1\leq p<+\infty$, and $u:\R^n\to\R$ be a measurable function that is $2\pi$-periodic in the variable $x_1$ and with finite seminorm $[u]_{W^{s,p}}^{\eper}.$ Then, 
\begin{equation}\label{ineq:scharz.rea}
  \left[u^{*\,n,1}\right]_{W^{s,p}}^{\eper}\leq [u]_{W^{s,p}}^{\eper}.
\end{equation}
Moreover, the following holds:
\begin{itemize}
\item [$(i)$] If $p>1$ then the equality in \eqref{ineq:scharz.rea} holds if and only if $u=\pm|u|$ and there exists $a\in\R^{n-1}$ such that $u(x_1,x')=\pm u^{{*\,n,1}}(x_1,x'+a)$  for a.e.\! 
$(x_1,x')\in\R\times\R^{n-1}$.
\item[$(ii)$] If $p=1$ and $u$ is continuous then the equality in \eqref{ineq:scharz.rea} holds if and only if $u=\pm|u|$ and for almost every $\tau> 0$ there exists $a_{\tau}\in\R^{n-1}$ such that the superlevel sets satisfy $$\{|u|>\tau\}=\{|u|>\tau\}^{{*\,n,1}}+(0,a_{\tau}),$$
up to redefining $u$ on a set of vanishing $n$-dimensional Lebesgue measure.
\end{itemize}
\end{theorem}

Let us start with the proof of   \Cref{thm:gy:periodic nonlocal Polya}, stating that the periodic Gagliardo seminorm decreases under periodic rearrangement. For this purpose it is convenient to write the periodic Gagliardo seminorm in a  different form using the Laplace transform. We shall denote  the Gamma function by
$\Gamma(\lambda):=\int_0^{+\infty}dt\, t^{\lambda-1}e^{-t}$ for $\lambda>0.$ By the change of variables $t\mapsto zt$ we obtain the Laplace transform of  the function $t\mapsto t^{\lambda-1}/\Gamma(\lambda),$ which is
$$
  z^{-\lambda}=\frac{1}{\Gamma(\lambda)}\int_0^{+\infty}dt\,t^{\lambda-1}e^{-zt}\quad\text{ for }z>0.
$$
Using this identity with 
$$
  \lambda=\frac{n+sp}{2}
  \quad\text{ and }\quad
  z=|x'-y'|^2+(x_1-y_1+2k\pi)^2,
$$
combined with the identity
$$
  \left([u]_{W^{s,p}}^{\eper}\right)^p
  =
  \int_{\R^{n-1}}dx'\int_{-\pi}^\pi dx_1\int_{\R^{n-1}}dy'\int_{-\pi}^\pi
  dy_1\sum_{k\in\mathbb{Z}}\frac{|u(x)-u(y)|^p}{|x-y+2k\pi e_1|^{n+sp}}
$$
for functions $u$ which are $2 \pi$-periodic in the variable $x_1$,
gives the following representation of the periodic Gagliardo seminorm:
\begin{equation}
 \label{eq:per Gagl seminorm by Phi}
   \left([u]_{W^{s,p}}^{\eper}\right)^p=\frac{1}{\Gamma(\lambda)}
  \int_0^{+\infty}dt\,t^{\lambda-1}
  \int_{\R^{n-1}}dx'
  \int_{\R^{n-1}}dy'\,
  \Phi(u,x',y',t)e^{-|x'-y'|^2t},
\end{equation}
where
$$
  \Phi(u,x',y',t):=\int_{-\pi}^{\pi} dx_1\int_{-\pi}^{\pi} dy_1\,
  |u(x_1,x')-u(y_1,y')|^p
  \sum_{k\in\Z}e^{-(x_1-y_1+2k\pi)^2t}.
$$

In view of representation \eqref{eq:per Gagl seminorm by Phi} it will be sufficient to show that the functional $\Phi$ does not increase if we replace $u$ by $u^{\ast\eper},$ for every $x',y'\in\R^{n-1}$ and every $t>0.$ Although this rearrangement inequality for $\Phi$ is now a $1$-dimensional problem, its proof relies on two nontrivial ingredients. The first one is \Cref{theorem:gy:non expansivity under J} which we will apply, for every given $t>0$,  
  to the kernel 
\begin{equation}
 \label{eq:def:g sum k def}
  z\mapsto g(z,t):=\sum_{k\in\mathbb{Z}} e^{-(z+2k\pi)^2t}.
\end{equation}
This will show that $\Phi$ does not increase under periodic rearrangement. Regarding the assumptions in \Cref{theorem:gy:non expansivity under J}, it is obvious that $g$ is $2\pi$-periodic and even with respect to the variable~$z$. The second nontrivial ingredient is the other hypothesis that $g$ must satisfy, to be decreasing in $(0,\pi)$ with respect to the variable $z$. This property follows from the monotonicity of the heat kernel on the circle, which in turns follows from a maximum principle on its spatial derivative based on the periodicity of the circle; see for instance \cite[Appendix B]{CCM Delaunay}.

\begin{remark}\label{remark:comparison n is one}{\em
In \cite[Theorem 1.4]{CCM Semilinear Var} we considered, for $n=1$ and $p=2$, a more general periodic Gagliardo seminorm where $|x-y|^{-(1+2s)}$ in \eqref{eq:def of periodic Gagl seminorm} is replaced by a more general kernel $K(x-y).$ However, the assumptions on $K$ are essentially such that the Riesz rearrangement inequality on the circle, i.e.,  \Cref{thm:gy:sharp Friedberg Luttinger}, can be directly applied.
For instance, a convexity assumption on $K$ gives that $g(t)=\sum_{k\in \Z}K(|t+2k\pi|)$ satisfies the assumption in  \Cref{thm:gy:sharp Friedberg Luttinger} on $g$ being symmetric and decreasing in $(0,\pi),$ and hence $g=g^{*\eper}.$ The second possible assumption on $K$ is such that after a Laplace transform one obtains a $g$ as in \eqref{eq:def:g sum k def} but with an additional nonnegative factor that only depends on $t.$ Obviously such kind of generalizations are also possible in the $n$-dimensional case.}
\end{remark}

\begin{proof}[Proof of  \Cref{thm:gy:periodic nonlocal Polya}.]
The proof will be divided into several steps.

\medskip
\noindent\underline{{\em Step 1} (proof of the inequality):} 
\smallskip

For the proof of \eqref{eq:periodic nonlocal Polya ineq}, we can assume that $u\geq 0$ due to the inequality $|u(x)-u(y)|\geq||u(x)|-|u(y)||$ and the fact that 
$u^{*\eper}=|u|^{*\eper}$, since $u^{*\eper}$ is defined by rearranging the superlevel sets of $|u|$.

In view of \eqref{eq:per Gagl seminorm by Phi}, it is sufficient to show that $\Phi$ does not increase under periodic rearrangement, for every $x',\,y'\in\R^{n-1}$ and every $t>0.$ If $g$ is defined as in \eqref{eq:def:g sum k def}, then we obtain from  \Cref{theorem:gy:non expansivity under J} 
---recall that $g(\cdot,t)=g(\cdot,t)^{*\eper}$ is decreasing in $(0,\pi)$; see  \cite[Appendix B]{CCM Delaunay}--- that
\begin{equation}
\begin{split}
  \Phi(u,x',y',t)&=\int_{-\pi}^{\pi}dx_1\int_{-\pi}^{\pi}dy_1\,
  |u(x_1,x')-u(y_1,y')|^pg(x_1-y_1,t)
   \\
  &\geq
  \int_{-\pi}^{\pi}dx_1\int_{-\pi}^{\pi}dy_1\,|u(\cdot,x')^*(x_1)-u(\cdot,y')^*(y_1)|^pg(x_1-y_1,t),
\end{split}
\end{equation}
where here $u(\cdot,x')^*$ denotes the Steiner symmetrization of $u(\cdot,x')\chi_{(-\pi,\pi)}$.
Now, using that for every $x'\in \R^{n-1}$ we have $u(\cdot,x')^*=u(\cdot,x')^{*\eper}$ in $(-\pi,\pi)$, and the identity \eqref{eq:astper reduce to 1 dim}, we obtain
\begin{equation}
 \label{eq:gy:Phi u x'y' t smaller}
  \Phi(u,x',y',t)\geq \Phi(u^{*\eper},x',y',t),
\end{equation}
which concludes the proof of the inequality.

\medskip
\noindent\underline{{\em Step 2} (case of equality for  $p>1$):} 
\smallskip

Firstly, observe that in case of equality $u$ cannot change sign, due to the following two facts: $|u|^{*\eper}=u^{*\eper}$ and $|u(x)-u(y)|>||u(x)|-|u(y)||$ on a set of positive measure in 
$\R^{n}\times\R^n$ if $u$ changes sign. Hence $u=\pm|u|.$ We assume from now on that $u\geq 0.$

In case of equality we must have that at least for one $t>0$ there is equality in \eqref{eq:gy:Phi u x'y' t smaller} ---actually there must be equality for almost every $t>0,$ but one $t$ is enough to conclude. Let us say for simplicity that $t=1$ and abbreviate $g(z)\equiv g(z,1).$  Thus, we must have that
\begin{equation}
\begin{split}
 \label{eq:gy:int LL u x1 xprime equals star}
  \intp dx_1\intp &dy_1\, |u(x_1,x')-u(y_1,y')|^pg(x_1-y_1)
   \\
  &=
  \intp dx_1\intp dy_1\, |u^{*\eper}(x_1,x')-u^{*\eper}(y_1,y')|^pg(x_1-y_1)
\end{split}
\end{equation}
for all $x'\in\R^{n-1}\setminus {N}$ and all $y'\in\R^{n-1}\setminus {N}({x'})$, where 
$\mathcal{L}^{n-1}({N})=\mathcal{L}^{n-1}({N}({x'}))=0$, and where $g$ satisfies the hypothesis of  \Cref{theorem:gy:non expansivity under J} and, in addition, $g$ is  decreasing in $(0,\pi)$. Let us denote
$$
  \Lambda(u):=\left\{x'\in\R^{n-1}:\, u(\cdot,x')=c(x')\in\R\text{ is constant in the $x_1$ variable}\right\}.
$$
Then, by  \Cref{theorem:gy:non expansivity under J} and \eqref{eq:astper reduce to 1 dim}, if $(x',y')$ is a pair satisfying the equality \eqref{eq:gy:int LL u x1 xprime equals star} and $x',y'\notin\Lambda(u),$ there exists $a(x',y')\in\R$ such that
\begin{align}
  u(x_1,x')&=u^{*\eper}(x_1+a(x',y'),x')\quad\text{ for a.e.\! }x_1\in\R,
   \\
  u(y_1,y')&=u^{*\eper}(y_1+a(x',y'),y')\quad\text{ for a.e.\! }y_1\in\R.
\end{align}
From this we will now deduce, since $x'\notin \Lambda(u),$ that $a(x',y')$ is independent of $x'$ and $y'.$ To see this, one can argue as follows. First of all, if $x'\notin \Lambda(u)$ then $u^{*\eper}(\cdot,x')$ is also nonconstant in $x_1$ and, therefore, $u^{*\eper}(x_1+b,x')=u^{*\eper}(x_1+d,x')$ forces that $b=d$ modulo $2\pi$; see the lines immediately after \eqref{eq:uS.uSper.tau} for a similar argument. Next, take another $z'\in\R^{n-1}\setminus\Lambda(u)$ different from $y'$ for which it also holds that $u(x_1,x')=u^{*\eper}(x_1+a(x',z'),x')$ for almost every $x_1\in\R.$ Then, $a(x',y')=a(x',z')$ modulo $2\pi$ and, therefore, $a$ is independent of the second variable. By symmetry, $a$ is also independent of the first variable. Thus, we have shown that $u(x_1,x')=u^{*\eper}(x_1+a,x')$ for almost every $x'\in\R^{n-1}\setminus \Lambda(u).$ But this also holds trivially for $x'\in\Lambda(u),$ and we obtain
$$
  u(x_1,x')=u^{*\eper}(x_1+a,x')\quad\text{ for a.e.\! }x'\in\R^{n-1}\text{ and a.e.\! }x_1\in\R.
$$

\medskip
\noindent\underline{{\em Step 3} (case of equality for  $p=1$):} 
\smallskip

Let us first address the case $(ii)$; the statement $(i)$ of \Cref{thm:gy:periodic nonlocal Polya} which deals with $p=1$ will be commented afterwards.

In case of equality, at the beginning of Step 2 we have seen that \eqref{eq:gy:int LL u x1 xprime equals star} holds for all $x'\in\R^{n-1}\setminus {N}$ and all $y'\in\R^{n-1}\setminus {N}({x'})$, where 
$\mathcal{L}^{n-1}({N})=\mathcal{L}^{n-1}({N}({x'}))=0$. However, since u is continuous, it is uniformly continuous on compact sets. The same holds true for $u^{*\eper}$, since Steiner symmetrization does not increase the modulus of continuity of a function; see \cite[Theorem 3.3, Section 4, and Remark 6.1]{Brock Solynin}. From this and a limiting argument it follows that if $u$ is continuous then \eqref{eq:gy:int LL u x1 xprime equals star} with $p=1$ actually holds for \textit{all} $x',y'\in\R^{n-1}$, and not only up to sets of measure zero.

As for $p>1,$ we can assume that $u\geq 0.$ Then, we obtain from  \eqref{eq:gy:int LL u x1 xprime equals star} and  \Cref{theorem:gy:non expansivity under J} that for every
\begin{equation}\label{def:RR(x,y)}
  \tau\in (\operatorname{ess}\inf u(\cdot,x'),\operatorname{ess}\sup u(\cdot,x'))\cap (\operatorname{ess}\inf u(\cdot,y'),\operatorname{ess}\sup u(\cdot,y'))=:\mathcal{R}(x',y')
\end{equation}
there exists  $a_{\tau}(x',y')\in \R$ such that
\begin{equation}
 \label{eq:gy:u( cdot xprime) bigger tau eq ast}
  \begin{split}
  \{u(\cdot,x')>\tau\}&= \{u^{*\eper}(\cdot,x')>\tau\}+a_{\tau}(x',y'),
 \\
  \{u(\cdot,y')>\tau\}&= \{u^{*\eper}(\cdot,y')>\tau\}+a_{\tau}(x',y')
 \end{split}
\end{equation}
in measure.
Consider for a moment only the first equation in \eqref{eq:gy:u( cdot xprime) bigger tau eq ast}. Note that  $a_{\tau}$  cannot depend on $y'$ because if 
$\tau\in (\essinf u(\cdot,x'),\esssup u(\cdot,x'))$ then 
$\{u^{*\eper}(\cdot,x')>\tau\}$ is a nontrivial ``periodic interval''  ---and hence all its translations which are not multiples of $2\pi$ are different. Thus, we obtain
$$
  \{u(\cdot,x')>\tau\}= \{u^{*\eper}(\cdot,x')>\tau\}+a_{\tau}(x')\quad\text{ for all $\tau\in\mathcal{R}(x',y')$}
$$
(this equality holds only in measure). 
In particular, and since \eqref{eq:gy:int LL u x1 xprime equals star} with $p=1$ holds for all $x',y'\in\R^{n-1}$, this last identity holds taking $x'=y'$, that is, for all $\tau\in \mathcal{R}(x',x')$ ---this is the key point where the continuity assumption on $u$ is used. Thus,
\begin{equation}
 \label{eq:gy:set u cdot xprime for all tau}
  \{u(\cdot,x')>\tau\}=\{u^{*\eper}(\cdot,x')>\tau\}+a_{\tau}(x')\quad\text{ for all }\tau>0.
\end{equation}
Next, we return to \eqref{eq:gy:u( cdot xprime) bigger tau eq ast}, which gives that 
$$
  \{u(\cdot,y')>\tau\}=\{u^{*\eper}(\cdot,y')>\tau\}+a_{\tau}(x')\quad\text{ for all }\tau\in \mathcal{R}(x',y').
$$
Thus $a_{\tau}$ does not depend on $x'.$ Now the claim follows immediately from the fact that $\{u>\tau\}$ is the union over all $x'\in\R^{n-1}$ of the sets $\{u(\cdot,x')>\tau\}\times\{x'\}$ (and from Fubini's theorem to deal with measure zero sets). In more detail, let $N(x')\subset \R$ be a set of zero $\mathcal{L}^1$-measure modulo which \eqref{eq:gy:set u cdot xprime for all tau} holds. Then, the set
$$
  M=\bigcup_{x'\in\R^{n-1}}N(x')\times\{x'\}
$$
has $\mathcal{L}^{n}$-measure zero by Fubini's theorem. Now by \eqref{eq:gy:set u cdot xprime for all tau}
$$
  \chi_{\{u(\cdot,x')>\tau\}}(x_1)= 
  \chi_{\{u^{*\eper}(\cdot,x')>\tau\}}(x_1-a_{\tau})\quad
  \text{ for all }x_1\in \R\setminus N(x').
$$
Hence this last identity holds for all $x\in\R^n\setminus M,$ which yields
$$
  \chi_{\{u>\tau\}}(x)=\chi_{\{u^{*\eper}>\tau\}}(x_1-a_{\tau},x')
$$
for all $x\in\R^n\setminus M$.

Finally, let us address the statement $(i)$ of \Cref{thm:gy:periodic nonlocal Polya} which deals with $p=1$. Recall that in case of equality, at the beginning of Step 2 we have seen that \eqref{eq:gy:int LL u x1 xprime equals star} holds for all $x'\in\R^{n-1}\setminus {N}$ and all $y'\in\R^{n-1}\setminus {N}({x'})$, where 
$\mathcal{L}^{n-1}({N})=\mathcal{L}^{n-1}({N}({x'}))=0$. 
Now, if $u$ is the characteristic function of a set $E$, we see that 
either $\mathcal{R}(x',y')=(0,1)$ or $\mathcal{R}(x',y')=\emptyset$ in \eqref{def:RR(x,y)}. Indeed, if $E$ has positive measure, for every $x'\in\R^{n-1}\setminus {N}$ we can always choose $y'\in\R^{n-1}\setminus {N}({x'})$ such that 
 $\mathcal{R}(x',y')=\mathcal{R}(x',x')$. Therefore, the rest of the proof follows as in \eqref{eq:gy:set u cdot xprime for all tau} and thereafter, showing that $u$ is of the form described in \Cref{thm:gy:periodic nonlocal Polya} $(ii)$. However, since $u$ is the characteristic function of a set, there is only one nontrivial superlevel set and, thus, $u$ is indeed of the form described in \Cref{thm:gy:periodic nonlocal Polya} $(i)$.
\end{proof}

We now prove the second main theorem of this section, which concerns the periodic Gagliardo seminorm and the cylindrical rearrangement.

\begin{proof}[Proof of  \Cref{thm:gy:Polya for n minus 1 Schwarz}.] 
We proceed as in the proof of  \Cref{thm:gy:periodic nonlocal Polya} and obtain using the Laplace transform as in \eqref{eq:per Gagl seminorm by Phi} the following representation of the periodic Gagliardo seminorm
$$
  ([u]_{W^{s,p}}^{\eper})^p
  =
  \frac{1}{\Gamma(\lambda)}\int_0^{+\infty}dt\,t^{\lambda-1}\int_{-\pi}^{\pi}dx_1 \int_{-\pi}^{\pi}dy_1\,\Psi(u,x_1,y_1,t)\sum_{k\in\Z}e^{-(x_1-y_1+2k\pi)^2 t},
$$
where
$$
  \Psi(u,x_1,y_1,t):=\int_{\R^{n-1}}dx'\int_{\R^{n-1}}dy'\,|u(x_1,x')-u(y_1,y')|^p e^{-|x'-y'|^2t}.
$$
By  \Cref{theorem:gy:J in R n-1} it holds that $\Psi(u^{{*\,n,1}},x_1,y_1,t)\leq \Psi(u,x_1,y_1,t)$ for almost every $x_1,y_1\in (-\pi,\pi)$ 
 and every $t>0$ (not necessarily for every $x_1$ and $y_1,$ since $f\in L^p((-\pi,\pi)\times\R^{n-1})$ only yields that $\int_{\R^{n-1}}dx'\,|f(x_1,x')|^p<+\infty$ for almost every $x_1$).
From this the inequality in the theorem follows.

The study of the case of equality is almost identical to that of  \Cref{thm:gy:periodic nonlocal Polya} and we will only carry out some details which are different.  If there is equality in the inequality then
\begin{equation}
\label{eq:gy:Psi u x1 y1 t}
  \Psi(u,x_1,y_1,t)=\Psi(u^{{*\,n,1}},x_1,y_1,t) 
\end{equation}
for all $x_1\in(-\pi,\pi)\setminus {N}$ and all $y_1\in(-\pi,\pi)\setminus {N}({x_1})$, where 
$\mathcal{L}^{1}({N})=\mathcal{L}^{1}({N}({x_1}))=0$), for some $t\in(0,+\infty)$ (we will need only one $t\in(0,+\infty)$ to conclude). 

Let us first address the case $p>1$. In this case we can argue identically as in the proof of  \Cref{thm:gy:periodic nonlocal Polya}, with the only difference that we replace $\Lambda(u)$ by 
$$
  \Lambda(u):=\{x_1\in (-\pi,\pi):\,u(x_1,x')=0 \text{ for a.e.\! $x'\in \R^{n-1}$}\}.
$$

We now address the case $p=1$. Using the continuity of $u$ ---as we did at the beginning of Step~3 in the proof of \Cref{thm:gy:periodic nonlocal Polya}--- together with \Cref{theorem:gy:J in R n-1}, we obtain that  for every $x_1,y_1\in (-\pi,\pi)$ and for all 
$$
  \tau\in \mathcal{R}(x_1,y_1):=(0,\min\{\operatorname{ess}\sup u(x_1,\cdot),
  \operatorname{ess}\sup u(y_1,\cdot)\})
$$
there exists $a_{\tau}(x_1,y_1)\in \R^{n-1}$ such that
\begin{align}
  \{u(x_1,\cdot)>\tau\}&=\{u^{{*\,n,1}}(x_1,\cdot)>\tau\}-a_{\tau}(x_1,y_1),
  \\
  \{u(y_1,\cdot)>\tau\}&=\{u^{{*\,n,1}}(y_1,\cdot)>\tau\}-a_{\tau}(x_1,y_1).
\end{align}
This two identities hold up to sets of $\mathcal{L}^{n-1}$ measure zero (which might depend on $x_1$ respectively $y_1$).
We can now argue similarly as in the proof of the case of equality for $p=1$ of \Cref{thm:gy:periodic nonlocal Polya} to conclude that $a_{\tau}(x_1,y_1)$ is independent of  $x_1$ and $y_1$. Finally, use that 
$$
  \{u>\tau\}=\bigcup_{x_1\in \R}\{x_1\}\times\{u(x_1,\cdot)>\tau\}
$$
and Fubini theorem (to deal with sets of measure zero) to obtain the claim of the theorem.
\end{proof}

\end{document}